\documentclass[pdflatex,sn-mathphys-num]{sn-jnl}

\usepackage[utf8]{inputenc}%
\usepackage[T1]{fontenc}%
\usepackage{graphicx}%
\usepackage{multirow}%
\usepackage{amsmath,amssymb,amsfonts}%
\usepackage{amsthm}%
\usepackage{mathrsfs}%
\usepackage{xcolor}%
\usepackage{textcomp}%
\usepackage{manyfoot}%
\usepackage{booktabs}%
\usepackage{array}%
\usepackage{algorithm}%
\usepackage{algorithmicx}%

\def\R{\mathbb{R}}

\def\Bu{\mathbb{B}}
\def\E{\mathbb{E}}

\def\C{\mathcal{C}}
\def\F{\mathcal{F}}
\def\H{\mathcal{H}}
\def\L{\mathcal{L}}
\def\M{\mathcal{M}}
\def\N{\mathcal{N}}

\def\SS{\mathrm{\mathcal{S}}}
\def\dist{\mathrm{dist}}

\def\z{\mathit{\hat{z}}}
\def\x{\mathit{\hat{x}}}
\def\v{\mathit{\hat{v}}}
\def\w{\mathit{\hat{w}}}
\def\dom{\mathrm{dom}\,}
\def\gph{\mathrm{gph}\,}

\newcommand{\NN}{\mathcal{N}}
\newcommand{\MM}{\mathcal{M}}
\newcommand{\sym}{\mathrm{sym}}
\newcommand{\skw}{\mathrm{skew}}
\newcommand{\norm}[1]{\lVert #1\rVert}

\theoremstyle{thmstyleone}%
\newtheorem{theorem}{Theorem}%
\newtheorem{corollary}[theorem]{Corollary}%

\theoremstyle{thmstyletwo}%
\newtheorem{remark}{Remark}%

\theoremstyle{thmstylethree}%

\begin{document}

\title[Asymmetric inexact progressive decoupling]{An asymmetric inexact progressive decoupling algorithm for linkage monotone inclusions with applications to stochastic variational inequalities}

\author[1]{\fnm{E. L.} \sur{Buscaglia}}\email{buscaglia@cifasis-conicet.gov.ar}

\author[2]{\fnm{P. A.} \sur{Lotito}}\email{plotito@exa.unicen.edu.ar}

\author*[1,3]{\fnm{L. A.} \sur{Parente}}\email{parente@cifasis-conicet.gov.ar}

\affil[1]{\orgname{CIFASIS (CONICET-UNR)}, \orgaddress{\street{Ocampo y Esmeralda}, \city{Rosario}, \postcode{2000}, \country{Argentina}}}

\affil[2]{\orgdiv{Facultad de Ciencias Exactas}, \orgname{Universidad Nacional del Centro de la Provincia de Buenos Aires}, \orgaddress{\street{Pinto 399}, \city{Tandil}, \postcode{7000}, \country{Argentina}}}

\affil[3]{\orgdiv{FCEIA}, \orgname{Universidad Nacional de Rosario}, \orgaddress{\street{Av. Pellegrini 250}, \city{Rosario}, \postcode{2000}, \country{Argentina}}}

\abstract{We introduce the Asymmetric Inexact Progressive Decoupling Algorithm
(AIPDA) for linkage monotone inclusions, combining Spingarn's partial inverse with
a possibly \emph{nonsymmetric} variable metric and a relative-error criterion based
on enlargements. AIPDA is the asymmetric inexact proximal point framework of
Buscaglia, Lotito, Parente and Solodov applied to the partial inverse, and inherits
its global convergence and local linear (superlinear) rates. It contains classical progressive
decoupling and progressive hedging and, under a blockwise error condition,
Rockafellar's inexact progressive decoupling algorithm from some iteration on. For
multistage stochastic variational inequalities, AIPDA yields the Asymmetric Inexact
Progressive Hedging Algorithm (AIPHA); for affine problems with polyhedral
constraints, a nonsymmetric metric makes each scenario subproblem solvable in closed
form and in parallel. On stochastic Cournot electricity markets, AIPHA is faster in
wall-clock time than tuned progressive hedging on the monotone instances tested.}



\maketitle

\section{Introduction}\label{sec:intro}

Monotone inclusions are a central topic in variational analysis, serving
as a framework that encompasses convex optimization, variational
inequalities, and saddle-point problems
\cite{RW:varan,HUL:01,fp:03,BC:17}. Given a Hilbert space
$\H$ with inner product $\langle\cdot,\cdot\rangle$ and corresponding norm $\|\cdot\|$, and a maximally monotone  set-valued operator $T:\H\rightrightarrows\H$, the
problem is
\begin{equation}\text{Find } x\in\H\;\text{ such that }\;\label{Inc1} 0\in T(x).\end{equation}
Here, monotonicity stands for $\langle v-w,z-x\rangle\geq 0$, for all $x,z\in\dom T$, $v\in  T(z)$ and $w\in T(x)$, while maximality means that the graph of $T$ is not properly contained in the graph of any other monotone operator, with $\gph T=\{(x,w):x\in\dom T,w\in T(x) \}$.
A classical tool for solving \eqref{Inc1} is the
Proximal Point Algorithm (PPA) \cite{mar:72,roc:76b}, which has served as a basis for developing a large number of approximation schemes (e.g., \cite{rock:76, BU:93, BGLS:06}) and iterates from the current iteration $z^k$ by defining
\begin{equation}\label{PPM}
z_{k+1}=R_k(z_k):=(I+c_kT)^{-1}(z_k),
\end{equation}
where $c_k$ is a regularization parameter and $R_k$ is the resolvent operator, which is single-valued and non-expansive due to the maximal monotonicity of T (see \cite{RW:varan}).

To exploit the possibilities of the scheme \eqref{PPM}, several ideas have been developed. A first proposal is \emph{to take advantage of the problem's structure}: specialize to optimization, variational inequalities, decomposition of sums, etc.
When $\H$ decomposes into two orthogonal subspaces, $\SS$ and $\SS^\bot$, a variation of \eqref{Inc1} is given by  the  \emph{ linkage problem}
\begin{equation}\label{linkage}
w\in T(x),\quad  \text{ with }x\in\SS \text{ and }w\in\SS^\bot. 
\end{equation}
Problem \eqref{linkage} can be cast in the form \eqref{Inc1} via the \emph{Spingarn's partial inverse} $T^\#$ and PPA leads to the Spingarn's method of partial inverses
\cite{spi:83}. In the finite-dimensional case, a closely related method is the  \emph{progressive decoupling algorithm}, introduced in \cite{Rock:19b} which, in its simplest form, starts with $x^0\in \SS$, $w^0\in \SS^\bot$ and iterates  in the form 
\begin{equation}\label{PDA}
\x^{k}=(T^k)^{-1}(0), \qquad x_{k+1}=P\x^k,\quad w_{k+1}=w_k-rP^\bot \x^k,
\end{equation}
where $T^k(x)=T(x)+r(x-x_k)-w_k,$ with  $P$ and  $P^\bot$ being the orthogonal projections onto $\SS$ and $\SS^\bot$, respectively. This scheme extends the progressive hedging algorithm of \cite{RW:91} (from which it takes its name), arising in stochastic programming and adapted in \cite{RS:18} to multistage stochastic variational inequalities \cite{RW:17}.

A second idea is \emph{inexactness}: solve subproblems approximately and handle errors. The first inexact proximal scheme was introduced in \cite{roc:76b} and required the error sequence to be summable for global convergence. Under strong regularity hypotheses on $T$ ($T^{-1}$ Lipschitz continuous at 0), local linear convergence rate was obtained if

\begin{equation}\label{error1}
\|z_{k+1}-R_k(z_k)\|\leq e_k\|z_{k+1}-z_k\|, \text{ with } \sum e_k<\infty.
\end{equation}

In practice, $\{e_k\}$ is chosen a priori and, since $e_k\to 0$, the subproblems must be solved with increasing precision, which raises the computational load. A relative-error approach has been introduced in \cite{solsva:Teps,solsva:hybrid} and unified in \cite{solsva:unif}, allowing a non-vanishing tolerance parameter but introducing a correction step. Among the applications of \cite{solsva:unif} is the inexact progressive hedging algorithm of \cite{BLP:24}, which is an extension of \cite{RS:18}.

A third line of work has enriched the PPA with \emph{variable
metrics}, in which a symmetric positive
definite (SPD) matrix is introduced to generalize the proximal parameter and induce iteration-dependent norms. Early ideas on this topic were related to quasi-Newton methods (see \cite{bgls:95,bq:98,bq:00}). The first general variable metric scheme for maximally monotone operators
\cite{pls:08}, a finite-dimensional extension of \cite{solsva:unif}, underlies
the decomposition method of \cite{lps:09} and the applications of
\cite{plmrs:11,LPRS:13}.
In contrast, \cite{rock:23b} considers summable errors and leads to the inexact
progressive decoupling algorithm of \cite{rock:24}, where the proximal parameter
in \eqref{PDA} is replaced by an SPD matrix $J_k$ and tolerances of type
\eqref{error1} are stated in the induced norms; \cite{Rock:26} refines it for
optimization, with a deeper convergence analysis and a theoretical implementation
for nonconvex operators. Nonsymmetric proximal terms appear in \cite{he:09}, an
exact method with a correction step, extended in \cite{dong:12} to relative-error
inexactness.
The recent framework of \cite{BLPS:25} deals with both symmetric and nonsymmetric proximal terms,  unifies the methods of \cite{pls:08, dong:12,rock:23b} and provides for the first time relationships between hybrid relative-error and
summable-error conditions. From the current iterate $z_k$, the Asymmetric Inexact Proximal Point Method (AIPPM) of \cite{BLPS:25} chooses  $c_k\geq c>0$, positive definite matrices $M_k$ and $N_k$ (with $M_k$ symmetric), and finds
$\varepsilon_k>0$,  together with points $\z^k$ and $\v^k$ such that
 \begin{equation}  \label{hippm}
\left\{
\begin{array}{rcl}
\v^k \in M_k^{-1}T^{\varepsilon_k}(\z^k),&&
\delta^k = c_k\v^k+N_k(\z^k-z_k),\\
\|\delta^k\|_{M_k}^2+2c_k\varepsilon_k&\leq&
\sigma_k^2\left(\|c_k\v^k\|_{M_k}^2+
\|N_k(\z^k-z_k)\|_{M_k}^2\right),\\ 
\langle \delta^k,\z^k-z_k\rangle_{M_k} + c_k\varepsilon_k&\leq&\sigma^2_k\langle N_k(\z^k-z_k),\z^k-z_k\rangle_{M_k},
\end{array}\right.
\end{equation}
where $\sigma_k\in[0,\bar\sigma)$ with $\bar\sigma<1$, $T^\varepsilon$ stands for the $\varepsilon$-enlargement of $T$ (see \cite{bis:97} and section 2 below)
and for a SPD matrix $M$ we denote $\langle
x,y\rangle_M =\langle Mx,y\rangle$ and $\|x\|_M = \sqrt{\langle Mx,x\rangle}$.
The update is defined by
\begin{equation}\label{Update1}z_{k+1} = z_k-\tau_ka_k \v^k, \quad a_k=\frac {\langle
\v^k,z_k-\z^k\rangle_{M_k}- \varepsilon_k} {\|
\v^k\|^2_{M_k}},\quad\tau_k\in(0,2).\end{equation}
Note that if $M_k=N_k=I$, $\sigma_k =0$ and $\tau_k=1$,  \eqref{PPM} is recovered. 
If $T^{-1}(0)\neq\emptyset$, the algorithm converges to an element
 $z^*\in T^{-1}(0)$  (see \cite[Theorem 4.3]{BLPS:25}) if:
 \begin{eqnarray}\label{Mcond0}
&&0<\lambda_M^l\leq\lambda_{min}(M_k)\leq \lambda_{max}(M_k) \leq \lambda_M^u,\\ \label{Mcond}
 &&M_{k+1}\preceq (1+\zeta_k) M_k,\quad \zeta_k >0 \,,  
 \quad \sum_{k=0}^{\infty} \zeta_k<\infty ,\\
\label{Ncond}
&&0<\sigma_N^l\leq\sigma_{min}(N_k)\leq \sigma_{max}(N_k)\le\sigma_N^u,\\&& M_kN_k \text{ is positive definite with }
  \sigma_N^l\leq \frac{1}{2}\lambda_{min}(N_k^TM_k+M_kN_k),\label{MNcond}
\end{eqnarray}
where for SPD matrices $A,B$, the notation $A\preceq B$ means that
$B-A$ is positive semidefinite, $\lambda_{min}(A)$ and $\lambda_{max}(A)$ are the minimal and maximal eigenvalues of $A$,  and for a PD matrix $N$, $\sigma_{min}(N)$ and $\sigma_{max}(N)$ are the minimal and maximal singular values of $N$.
If
 $T^{-1}$ is calm at zero for $z^*$, the linear (superlinear) rate of convergence is obtained under appropriate choices of the involved 
parameters (see \cite[Theorem 4.4]{BLPS:25}) or by strengthening \eqref{Mcond} to 
 \begin{equation}\label{Mcond1}
\frac{1}{1+\zeta_k}M_k\preceq M_{k+1}\preceq(1+\zeta_k)M_k,\quad \zeta_k >0 \,,
\quad \sum_{k=0}^{\infty} \zeta_k<\infty .
\end{equation}

The present article combines these three ideas, drawing on \cite{rock:24} and
\cite{BLPS:25}. We introduce the \emph{Asymmetric Inexact Progressive Decoupling
Algorithm} (AIPDA, Algorithm~\ref{AIPDA}), which couples the partial inverse of
\cite{spi:83} with symmetric and nonsymmetric variable metrics and a relative-error
criterion that uses the enlargement of $T$, and we prove that its iterates coincide
with those of the AIPPM of \cite{BLPS:25} applied to $T^{\#}$. Hence AIPDA inherits
the global convergence and local linear (superlinear) rate of \cite{BLPS:25}. With
scalar metrics AIPDA reduces to the progressive decoupling algorithm of
\cite{Rock:19b} and, in the stochastic setting, to progressive hedging
\cite{RW:91,RS:18}. The inexact progressive decoupling algorithm of \cite{rock:24}
(and its optimization version \cite{Rock:26}) uses symmetric metrics and summable errors; we show that, under a blockwise version of
its error condition, its iterations are AIPDA iterations from some index on. AIPDA
thus adds nonsymmetric metrics, a free parameter $c_k$ (which allows superlinear
convergence) and relative errors with a non-vanishing tolerance. Finally, specializing
AIPDA to multistage stochastic variational inequalities yields the \emph{Asymmetric
Inexact Progressive Hedging Algorithm} (AIPHA, Algorithm~\ref{AIPHA-eps}), which
extends \cite{BLP:24} to nonsymmetric variable metrics and which we compare with
progressive hedging on stochastic Cournot electricity markets.

Section~\ref{sec:prelim} collects background on partial inverses, enlargements and
regularity conditions, Section~\ref{sec:algorithm} presents AIPDA and its equivalence
with AIPPM, Section~\ref{sec:aipha} develops AIPHA, Section~\ref{sec:numerics}
reports the numerical experiments and Section~\ref{sec:conclusions} concludes.

\section{Preliminaries}\label{sec:prelim}

Let $\H$ be a finite-dimensional Hilbert space, with $\SS$ being a
nonempty subspace. Define $P$ as the orthogonal projection onto $\SS$ and
$P^\bot$ as the orthogonal projection onto $\SS^\bot$, the orthogonal
complement of $\SS$.
For a maximally monotone operator $T:\H\rightrightarrows\H$, the partial inverse of $T$ with respect to $\SS$, denoted by $T^\#$, is given by
$$v\in T^\#(z)\Leftrightarrow Pv+P^\bot z\in T(P z + P^\bot v).$$
Let $J:\H\to\H$ be a linear operator preserving  $\SS$
and $\SS^\bot$, i.e., $JPu\in\SS$ and $JP^\bot u\in\SS^\bot$, for all
$u\in\H$.
Then, $PJ=JP$, $P^\bot J=JP^\bot$ and, if  $J$ is invertible, it holds that
\begin{equation}\label{PartMatrix} J^\# z =  J P z+ J^{-1}P^\bot z.\end{equation}
On the other hand,  for $\varepsilon\geq 0$, the $\varepsilon$-enlargement of $T$ given by
\begin{equation*}
    T^\varepsilon(z) = \{ v\in \H \, | \, \langle v - x, z-w\rangle \geq -\varepsilon,\ \forall w\in  \H, \forall x \in T(w)\}.
\end{equation*}
Several examples and properties of $\varepsilon$-enlargements can be found in \cite{bis:97}; in particular,  $T^0= T$ and if $\varepsilon=\varepsilon_1+\varepsilon_2$, both nonnegative, then
\begin{equation}(A^{\varepsilon_1}+B^{\varepsilon_2})\subset (A+B)^\varepsilon.\label{enlargesum}\end{equation} 
Also, it is well known that enlargements and partial inverses commute \cite[Lemma 3.1]{BSS:06}, i.e.,
\begin{equation}(T^\#)^\varepsilon=(T^\varepsilon)^\#.\label{partinven}\end{equation}

Local regularity conditions of Lipschitz type for set-valued operators are established in terms of its graph. We give here some basic definitions and facts, the reader may consult \cite{DR:09,Ioffe:16,KLT:17,KLT:18} for a comprehensive approach. 

A mapping $T: \H \rightrightarrows \H$ is said to be outer-Lipschitz-continuous at $\bar x \in \dom T $ if $T(\bar x)$ is a closed set, and there exist a constant $L\geq 0$ along with a neighborhood $U$ of $\bar x$ such that $T(x) \subset T(\bar x) +L\|x-\bar x\|\Bu$ for all $x \in U $, where $\Bu$ is the closed unit ball in $\H$.

A mapping $T: \H \rightrightarrows \H$ is \emph{metrically subregular} at $x^*$ for $y^*$ if 
$(x^*, y^*) \in \gph\ T$ 
and there exist neighborhoods $U$ of $x^*$ and $V$ of $y^*$, along with $L>0$, such that for all $x\in U$ it holds that $ \dist (x, T^{-1}(y^*)) \leq L \, \dist (y^*, T(x)\cap V).$

A mapping $S: \H \rightrightarrows \H$ 
is said to be \emph{calm} at $y^*$ for $x^*$ if $(y^*, x^*) \in \gph S$, and
there exist a constant $L\geq 0$ along with neighborhoods $V$ of $y^*$ and $U$ of $x^*$ such that 
$S(y)\cap U \subset S(y^*) + L \|y-y^*\|\Bu$ for all $y \in V$.

It is clear from the definitions that Lipschitz continuity implies calmness. Also, $T$ is metrically subregular at $ x^*$ for $y^*$ if and only if $T^{-1}$ is calm at
 $y^*$ for $x^*$ (see \cite[Theorem 3H.3]{DR:09}).

A closely related notion is subtransversality.  Closed sets $A,B$ in $\H$ are subtransversal at $x^*\in A\cap B$  if there is a constant $L>0$ and a neighborhood $U$ of $x^*$ such that for any $x\in U$, it holds that \footnote{In some texts, like \cite{Ioffe:16}, subtransversality is stated in the  form $d(x,A\cap B)\leq K \{d(x,A)+d(x,B)\}$, which is equivalent to \eqref{subtransv} for proper choices of $K$.}
\begin{equation}\label{subtransv}
d(x,A\cap B)\leq L\max \{d(x,A),d(x,B)\}. 
\end{equation}

Metric subregularity of an operator in $\H$ connects with subtransversality of sets in $\H\times \H$ as follows. Given a closed-graph set-valued mapping $T:\H\rightrightarrows\H$ and $(x^*,y^*)\in\gph T$, let $A=\gph T$ and $B=\H\times\{y^*\}$, so $(x^*,y^*)\in A\cap B\subset\H\times\H$.
Then $A,B$ are subtransversal at $(x^*,y^*)$ if and only if $T$ is metrically subregular at $x^*$ for $y^*$ (see \cite[Proposition 3]{KLT:17}). 

For finite-dimensional Hilbert spaces, subtransversality holds globally for convex polyhedral sets. A convex polyhedral set $Q$ is a finite intersection of closed half-spaces and hyperplanes.
If $Q_1$ and $Q_2$ are convex polyhedral, they are globally subtransversal, i.e., there exists $L>0$ such that
\begin{equation}\label{globsubt}
    d(x,Q_1\cap Q_2)\leq L\max\{d(x,Q_1),d(x,Q_2)\},
\end{equation}
(see \cite[Proposition 7.28]{Ioffe:16}).
A set-valued mapping $T:\H\rightrightarrows\H$ is called \emph{polyhedral} if $\gph T$ is the union of finitely many convex polyhedral sets of $\H\times\H$; then $T^{-1}$ is polyhedral as well. A classical result of Robinson \cite{rob:81} states that every polyhedral mapping $T$ is outer Lipschitz at every $\bar y$. In particular, if $T^{-1}$ is polyhedral, then $T^{-1}$ is calm at $y^*$ for every $x^*\in T^{-1}(y^*)$. Typical examples are $T=F+N_C$ with $F$ affine and $C$ convex polyhedral, since $\gph N_C$ is a finite union of convex polyhedral sets (one for each face of $C$).

\section{The AIPDA algorithm}\label{sec:algorithm}
 
 In this section, the new progressive decoupling algorithm is introduced, an inexact scheme that trivially extends the method in \cite{Rock:19b} (described in \eqref{PDA}). The framework of \cite{BLPS:25} gives the basis to analyze its convergence properties and to establish relationships with other progressive decoupling methods.
 
 \begin{algorithm}[H]
 \caption{Asymmetric Inexact Progressive Decoupling Algorithm (AIPDA)}\label{AIPDA}
\begin{algorithmic}[1]

\item \textit{Initialization:} Choose $x_0\in \SS$, $w_0\in \SS^\bot, \  \bar{\sigma} \in (0,1), \theta \in (0,1)$ and $c>0$.

\item \textit{Inexact Proximal Step:} Choose $\sigma_k\in [0,\bar{\sigma})$, $c_k\geq c$, a self-adjoint positive-definite linear mapping $S_k$ and a positive-definite linear  mapping $A_k$,  both preserving $\SS$ and $\SS^\bot$. 
Define $D_k=c_k^{-1}P+c_kP^\bot$.
Find $\hat x^k$, $\hat w^k$ and $\varepsilon_k\geq 0$ such that
\begin{equation}\label{ips} \left\lbrace \begin{array}{l}
D_k A_k(x_k-\x^k)+w_k \in T^{\varepsilon_k}  (\w^k),\\
\hat \delta^k=\w^k-\x^k,
\end{array}
\right.\end{equation}
satisfying the inequalities
\begin{equation} \begin{array}{l}
\|(A_kP+c_kP^\bot)\hat \delta^k\|_{S_k}^2 +2c_k\varepsilon_k  \leq  \\\qquad\; \quad\sigma_k^2\big(\| A_k(x_k-P\x^k)+c_kP^\bot \w^k\|_{S_k}^2+\| A_k(x_k-P\w^k)+c_kP^\bot \x^k\|_{S_k}^2\big)
\label{eq:conddeltaeps1}\end{array}
\end{equation}
and
\begin{equation} \begin{array}{l}
\langle \hat\delta^k,c_k^{-1}A_k^\top(P\w^k-x_k)-c_kA_k P^\bot\x^k\rangle+\varepsilon_k\leq \\\qquad\; 
\sigma_k^2\langle D_kA_k(P\w^k-x_k+P^\bot\x^k),P\w^k-x_k+P^\bot\x^k\rangle.
\label{eq:conddeltaeps2}\end{array}
\end{equation}
\item \textit{Update:} If $x_k-P\w^k=P^\bot \x^k=0$, STOP. Otherwise, choose  $\tau_k\in [1-\theta,1+\theta]$ and set
 \begin{equation}\label{updates}
 \begin{array}{lcl}
 x_{k+1}&=& x_k - \tau_k \alpha_k S_kA_k(x_k-P\x^k)
\\
 w_{k+1}&=&w_k-\tau_k \alpha_k c_k S_kP^\bot\w^k,
 \end{array}
 \end{equation}
where
\begin{equation*}
    \alpha_k
   = \frac{\langle A_k(x_k-P\x^k+c_kP^\bot\x^k),x_k -P\w^k +c_k P^\bot \w^k\rangle - c_k\varepsilon_k} {\|A_k(x_k-P\x^k)+c_kP^\bot \w^k\|_{S_k}^2}
\end{equation*}
Set $k:=k+1$ and go to the Inexact Proximal Step.
\end{algorithmic}
\end{algorithm}

\begin{remark}\label{rem:PDA}
Note that, if the proximal step is solved exactly in the previous algorithm, (that is,  $\sigma_k\equiv 0$) then \eqref{eq:conddeltaeps1} implies $\varepsilon_k\equiv 0$ and $\hat \delta_k\equiv 0$, so \eqref{eq:conddeltaeps2} is trivially fulfilled (for any choice of $S_k$) and \eqref{ips} is an exact proximal step with $\x^k=\w^k$, for all $k$. If, in addition, $c_k\equiv 1$, $S_k\equiv r^{-2}P+P^\bot$ and $A_k\equiv rI$, then $\alpha_k\equiv r$ so that, choosing $\tau_k\equiv 1$, we obtain 
$\tau_k \alpha_k c_k S_kA_kP\equiv P$ and $\tau_k \alpha_k c_k S_kP^\bot\equiv rP^\bot$.
Hence, the updates in \eqref{updates} give $x^{k+1}=P\x^k$  and $w^{k+1}=w_k-rP^\bot\x^k$, for all $k$, so  the exact progressive decoupling algorithm \eqref{PDA} is recovered. 
Multiplying $S_k$ by a positive constant only rescales $\alpha_k$ and leaves the iterates unchanged. 
\end{remark}

The following result places Algorithm \ref{AIPDA} within the framework of \cite{BLPS:25}.

\begin{theorem}\label{thm:equiv}
The iterations given by the AIPDA  are equivalent to those obtained by the AIPPM from \cite{BLPS:25}, described in \eqref{hippm}-\eqref{Update1}, applied to the mapping $ T^\#$, by setting the matrices $M_k= S_k^{-1}$ and $N_k=S_kA_k^\#$, with the same parameters $c_k,\sigma_k$ and $\tau_k$, for all $k$.
\end{theorem}

\begin{proof}
    Let $x_k\in\SS$ and $w_k\in\SS^\bot$ be the current iterates of the AIPDA, and define $z_k=x_k+w_k$.
Then $P(z_k)=x_k$ and $P^\bot(z_k)=w_k$.
Observe that if $0\in T^\#(z_k)$ then  $0\in T^\#(x_k+w_k)$ and the definition of $T^\#$ implies that $w_k\in T(x_k)$ so $(x_k,w_k)$ solves the linkage problem \eqref{linkage}; hence the problem reduces to finding a zero of $T^\#$.
Consider an iteration of the AIPPM given by \cite{BLPS:25},  with $M_k= S_k^{-1}$, $N_k=S_kA_k^\#$, the same proximal parameter $c_k$ and the same tolerance parameter $\sigma_k$, that is, find $\v^k,\z^k$, $\varepsilon_k$ such that
\begin{equation}\label{IPS1}
\v^k \in  S_k(T^\#)^{\varepsilon_k}(\z^k),\qquad
\delta^k  =c_k \v^k +S_kA_k^\#(\z^k-z_k), 
\end{equation}
verifying
\begin{eqnarray}
\|\delta^k\|_{S_k^{-1}}^2 + 2c_k\varepsilon_k &\leq&  \sigma_k^2  \big( \|c_k\v^k\|_{S_k^{-1}}^2 + \|S_kA_k^\#(\z^k- z_k)\|_{S_k^{-1}}^2\big),
 \label{eq:conddelprox1}\\
 \langle\delta^k,\z^k- z_k\rangle_{S_k^{-1}} +c_k\varepsilon_k&\leq&\sigma_k^2\langle S_kA_k^\#(\z^k- z_k),\z^k- z_k\rangle_{S_k^{-1}}\label{eq:conddelprox2}
.\end{eqnarray}
The algorithm stops if $z_k-\z^k=0$ (in which case $z_k$ is a zero of $T^\#$ and thus $(x_k,w_k)$ solves the monotone inclusion), otherwise, select $\tau_k \in [1-\theta,1+\theta]$ and define
\begin{equation}z_{k+1}= z_k-\tau_ka_k\v^k,\quad a_k= \frac{\langle \v^k,z_k-\z^k\rangle_{S_k^{-1}} - \varepsilon_k}{\|\v^k\|_{S_k^{-1}}^2}.\label{update1}
\end{equation}
From \eqref{IPS1} we have
\begin{equation*}
  c_k^{-1}S_k^{-1}(\delta^k +S_kA_k^\#(z_k-\z^k))\in (T^\#)^{\varepsilon_k}(\z^k),
\end{equation*}
which,  from \eqref{partinven} and the definition of $T^\#$, gives
$$\begin{array}{l}
 c_k^{-1}PS_k^{-1}(\delta^k +S_kA_k^\#(z_k-\z^k))+ P^\bot\z^k\in\\ \qquad T^{\varepsilon_k}\left(P\z^k+ c_k^{-1}P^\bot S_k^{-1}(\delta^k +S_kA_k^\#(z_k-\z^k))\right).
\end{array}$$
From the linearity, the properties of the projections and the fact that $S_k$ and $A_k$ preserve the subspaces $\SS$ and $\SS^\bot$, using \eqref{PartMatrix} and separating terms, the above inclusion turns into
\begin{equation}
     c_k^{-1}[S_k^{-1}P\delta^k + A_k(x_k - P\z^k)] + P^\bot \z^k \in T^{\varepsilon_k}\left(P\z^k +  c_k^{-1}[S_k^{-1}P^\bot \delta^k + A_k^{-1} (w_k- P^\bot\z^k)]\right).
    \label{expandedInclusion}
\end{equation}
Define $\w^k:= P\z^k +  c_k^{-1}S_k^{-1}P^\bot \delta^k +  c_k^{-1}A_k^{-1} (w_k - P^\bot\z^k)$. Then
\begin{equation}P\w^k = P\z^k   \label{zhProj}
\end{equation}
$$ P^\bot\w^k=   c_k^{-1}S_k^{-1}P^\bot \delta^k +  c_k^{-1}A_k^{-1} w_k -  c_k^{-1}A_k^{-1} P^\bot\z^k$$
so
\begin{equation}\label{Pozhat} P^\bot\z^k = A_kS_k^{-1}P^\bot\delta^k+ w_k - c_kA_k P^\bot\w^k.\end{equation}
Replacing in \eqref{expandedInclusion}, we obtain
 \begin{equation}
     c_k^{-1}S_k^{-1}P\delta^k + c_k^{-1} A_kx_k -  c_k^{-1} A_kP\w^k +A_kS_k^{-1}P^\bot\delta^k+ w_k - c_kA_k P^\bot\w^k\in T^{\varepsilon_k}(\w^k)
 \end{equation}
 which gives
  \begin{equation}
    D_kA_k(x_k + A_k^{-1} S_k^{-1}P\delta^k +  c_k^{-1}S_k^{-1}P^\bot\delta^k-\w^k) + w_k\in T^{\varepsilon_k}(\w^k)
 \end{equation}
Setting \begin{equation}\label{xhat}
    \x^k := - A_k^{-1} S_k^{-1}P\delta^k - c_k^{-1}S_k^{-1}P^\bot\delta^k+\w^k,
\end{equation} 
the previous inclusion coincides with the one in \eqref{ips},
and it holds that \begin{equation}\label{deltahat}
    \delta^k = S_kA_kP(\w^k-\x^k)+c_kS_kP^\bot (\w^k-\x^k)=S_k(A_kP+c_kP^\bot)\hat\delta^k.
\end{equation}
From \eqref{Pozhat} and \eqref{deltahat}, it holds that
\begin{eqnarray}  \z^k -z^k &=& P(\z^k -z^k)+P^\bot (\z^k -z^k)\nonumber \\
    &=& P\w^k-x_k + A_kS_k^{-1}P^\bot\delta^k+ w_k - c_kA_k P^\bot\w^k - w_k\nonumber\\
    &=& P\w^k-x_k - c_kA_k P^\bot\x^k. \label{zhatmz}
\end{eqnarray}
Note that, since $x_k\in\SS$ and $A_k$ is positive definite, $\z^k-z_k=0$ if and only if $P\w^k-x_k=0$ and $P^\bot\x^k=0$ so both stopping criteria coincide.
In addition, using \eqref{PartMatrix},
\begin{equation}S_kA_k^\#(\z^k -z^k) = S_k(A_k(P\w^k-x_k)-c_kP^\bot \x^k),\label{Nzhatmz}
\end{equation}
so
\begin{eqnarray}
        c_k\v^k &=& \delta^k + S_kA_k^\#(z^k -\z^k)\nonumber  \\
         &=& \delta^k -S_kA_k(P\w^k-x_k)+c_kS_kP^\bot \x^k\nonumber \\
          &=& S_kA_kP(\w^k-\x^k)+c_kS_kP^\bot (\w^k-\x^k)-S_kA_k(P\w^k-x_k)+c_kS_kP^\bot \x^k\nonumber \\
          &=&S_k(A_k(x_k-P\x^k)+c_kP^\bot \w^k).\label{vhat}
\end{eqnarray}
Replacing \eqref{deltahat}, \eqref{vhat} and \eqref{Nzhatmz} in \eqref{eq:conddelprox1} and \eqref{eq:conddelprox2}, it holds that $ \delta^k$ verifies the conditions
$$ \begin{array}{l}  \|S_k(A_kP+c_kP^\bot)\hat\delta^k\|_{S_k^{-1}}^2 + 2 c_k\varepsilon_k\leq \\ \qquad\qquad \sigma_k^2\big(\|S_k(A_k(x_k-P\x^k)+c_kP^\bot\w^k)\|_{S_k^{-1}}^2 +\|S_k(A_k(x_k-P\w^k)+c_kP^\bot \x^k)\|_{S_k^{-1}}^2 \big)\end{array}$$
and
$$  \begin{array}{l} \langle S_k(A_kP+c_kP^\bot)\hat\delta^k,P\w^k-x_k - c_kA_k P^\bot\x^k\rangle_{S_k^{-1}} +c_k\varepsilon_k\leq  \\ \qquad\qquad\sigma_k^2\langle S_k(A_k(P\w^k-x_k)-c_kP^\bot \x^k) ,P\w^k-x_k-c_kA_k P^\bot\x^k\rangle_{S_k^{-1}}. \end{array}$$
This gives
$$ \begin{array}{l}  \|(A_kP+c_kP^\bot)\hat\delta^k\|_{S_k}^2 + 2 c_k\varepsilon_k\leq \\ \qquad\qquad \sigma_k^2\big(\|A_k(x_k-P\x^k)+c_kP^\bot\w^k\|_{S_k}^2 +\|A_k(x_k-P\w^k)+c_kP^\bot \x^k\|_{S_k}^2 \big)\end{array}$$
and
$$  \begin{array}{l} \langle \hat\delta^k,A_k^\top(P\w^k-x_k) - c_k^2A_k P^\bot\x^k\rangle +c_k\varepsilon_k\leq  \\ \qquad\qquad\sigma_k^2(\langle A_k(P\w^k-x_k) ,P\w^k-x_k\rangle+c_k^2\langle P^\bot\x^k,A_k P^\bot\x^k\rangle). \end{array}$$
The first inequality is just \eqref{eq:conddeltaeps1}.  Dividing by $c_k$ and rearranging terms (by using the orthogonality), the second one turns into
$$  \begin{array}{l} \langle \hat\delta^k,c_k^{-1}A_k^\top(P\w^k-x_k) -c_kA_k P^\bot\x^k\rangle +\varepsilon_k\leq  \\ \qquad\qquad\sigma_k^2\langle (c_k^{-1}PA_k+c_kP^\bot A_k)(P\w^k-x_k+ P^\bot\x^k) ,P\w^k-x_k+P^\bot\x^k\rangle,
\end{array}$$
which is equivalent to \eqref{eq:conddeltaeps2}.
Regarding $\alpha_k$,
\begin{equation}
\begin{split}
a_k &= \frac{\langle \v^k,z_k-\z^k\rangle_{S_k^{-1}} - \varepsilon_k} {\|\v^k\|_{S_k^{-1}}^2}  \\
   &= \frac{\langle c_k^{-1}S_k(A_k(x_k-P\x^k)+c_kP^\bot \w^k),x_k -P\w^k+ c_kA_k P^\bot\x^k\rangle_{S_k^{-1}} - \varepsilon_k} {\|c_k^{-1}S_k(A_k(x_k-P\x^k)+c_kP^\bot \w^k)\|_{S_k^{-1}}^2} \\
   &= c_k\frac{\langle A_k(x_k-P\x^k),x_k -P\w^k\rangle +\langle c_kP^\bot \w^k,c_kA_k P^\bot\x^k\rangle - c_k\varepsilon_k} {\|A_k(x_k-P\x^k)+c_kP^\bot \w^k\|_{S_k}^2} \\
   &= c_k\frac{\langle A_k(x_k-P\x^k+c_kP^\bot\x^k),x_k -P\w^k +c_k P^\bot \w^k\rangle - c_k\varepsilon_k} {\|A_k(x_k-P\x^k)+c_kP^\bot \w^k\|_{S_k}^2}=c_k\alpha_k. 
\end{split}
\end{equation}

Then the update $z_{k+1}$ can be written as
\begin{equation}
\begin{split}
    z_{k+1}&= z_k-\tau_k a_k\v^k \\
    & = x_k + w_k -\tau_k \alpha_k S_k(A_k(x_k-P\x^k)+c_kP^\bot \w^k) \\
    & = x_k - \tau_k \alpha_k S_kA_k(x_k-P\x^k) + (w_k-\tau_k \alpha_kc_kS_kP^\bot \w^k)\\
    &= x_{k+1}+ w_{k+1},
\end{split}
\end{equation}
which are the same updates given in \eqref{updates}. Note that the preceding procedure establishes a chain of equivalences, so it is reversible: starting with the current iterate $z^k$ of AIPPM, defining $x_k=Pz^k$, $w_k=P^\bot z^k$ and performing an iteration of AIPDA, the iterates of AIPPM are recovered by defining $\z^k$ as in \eqref{zhProj}-\eqref{Pozhat} and $\v^k$ as in \eqref{vhat}. Therefore,  both algorithms coincide.
\end{proof}

Theorem \ref{thm:equiv} says that AIPDA is a particular case of AIPPM, hence, under appropriate assumptions, it inherits its convergence properties: global convergence of the iterates and local linear (superlinear) convergence rate. The linear operators $S_k,A_k$ must satisfy hypotheses equivalent to \eqref{Mcond0}-\eqref{Mcond1}, that is, there exist parameters $0<\lambda_l\leq\lambda_u$, $0<\sigma_l\leq\sigma_u$, $0<\nu_l$, and a summable series of positive terms $\sum \zeta_k$, such that, for all $k$,

\begin{eqnarray}
&& \lambda_l\le\lambda_{min}(S_k)\le\lambda_{max}(S_k)\le\lambda_u, \label{H0new}  \\
&& S_{k+1}\succeq(1+\zeta_k)^{-1}S_k, \label{Hnew} \\
&&\sigma_l\le\sigma_{\min}(A_k)\le\sigma_{\max}(A_k)\le\sigma_u,\ \text{and}\
  \lambda_{min}\big(\tfrac12(A_k+A_k^{\top})\big)\ge\nu_l>0, \label{H2new} \\
 && (1+\zeta_k)^{-1}S_k\preceq \ S_{k+1}\preceq(1+\zeta_k)S_k.\label{H3new}
\end{eqnarray}
On the other hand, in order to achieve local convergence rate results, we ask for a regularity condition on the set-valued mapping $T$ in relation with the subspaces $\SS$ and $\SS^\bot$. 
In the following, let $Z^*$ be the solution set of the linkage problem \eqref{linkage}, that is, 
$$Z^*:=\{(x,w): x\in\SS, w\in\SS^\bot \text{ and } w\in T(x)\}\subset \H\times\H.$$

\begin{corollary}\label{ConvProp}
Suppose that $Z^*\neq\emptyset$ and conditions \eqref{H0new}, \eqref{Hnew} and \eqref{H2new}
hold. Then,  following statements hold.
\begin{enumerate}
    \item Any sequence $\{(x^k,w^k)\}$ generated by Algorithm \ref{AIPDA} converges to some $(x^*,w^*)\in Z^*$.
    \item 
In addition, assume that $\gph T$  and $\SS\times\SS^\bot$  are subtransversal at $(x^*,w^*)$. Then,
\begin{enumerate}
    \item If for $k$ large enough, the linear operators $S_k$ are chosen sufficiently close to the
identity, then the sequence $\{(x^k,w^k)\}$ converges R-linearly.
\item If the products $S_kA_k^\#$ are sufficiently close to the identity matrix for $k$ large
enough, then $\sigma_k$ can be chosen small enough, $c_k$ can be chosen large enough
and $\tau_k$ can be chosen sufficiently close to 1, so that the sequence  $\{(x^k,w^k)\}$ converges R-linearly or even Q-linearly. If $S_kA_k^\#\to I$, $\sigma_k\to 0$, $c_k\to\infty$ and
$\tau_k\to 1$, then convergence is superlinear.
\item If condition \eqref{Hnew} is strengthened to \eqref{H3new}, there exists $k_0\in\NN$  such that the sequence $\{(x^k,w^k)\}$
converges R-linearly in the norm induced by $S_{k}$, for any $k\ge k_0$, and in particular in the norm induced by $S_\infty:=\lim S_k$.
\end{enumerate}
\item In the case that $T$ is polyhedral then the hypothesis of subtransversality in item~2 holds at every
$(x^*,w^*)\in Z^*$, so (a)-(c) hold for the limit point of any sequence generated
by Algorithm~\ref{AIPDA}, with no further regularity assumption on $T$.
\end{enumerate}
\end{corollary}
\begin{proof}
    Note that $(x^*,w^*)\in Z^*$ if and only if $0\in T^\#(x^*+w^*)$. In view of Theorem \ref{thm:equiv}, the first result holds by applying \cite[Theorem 4.3]{BLPS:25}, if we prove that conditions \eqref{Mcond0}-\eqref{MNcond} hold, while the rest cames from \cite[Theorem 4.4 and Remark 4.1]{BLPS:25} if $(T^\#)^{-1}$ is calm at 0 for $(x^*+w^*)$.

    Since $M_k=S_k^{-1}$, conditions \eqref{Mcond} and  \eqref{Hnew} are equivalent, and condition \eqref{H0new} implies  $\lambda_M^l:=1/\lambda_u\le\lambda_{min}(M_k)\le\lambda_{max}(M_k)\le 1/\lambda_l=:\lambda_M^u$, so \eqref{Mcond0} holds.
    
    On the other hand, $A_k$ preserves $\SS$ and $\SS^\bot$. Hence, it is block diagonal with respect to the orthogonal decomposition $\H=\SS\oplus\SS^\bot$, which means that $A_k=A_{k,\SS}P+A_{k,\SS^\bot}P^\bot$, and the singular values of $A_k$ are those of its two blocks $A_{k,\SS}$ and $A_{k,\SS^\bot}$. By \eqref{PartMatrix},
 $A_k^{\#}=A_{k,\SS}P+A_{k,\SS^\bot}^{-1}P^\bot$,
so, using \eqref{H2new}, the singular values of $A_k^{\#}$ lie in $[\min\{\sigma_l,1/\sigma_u\},\max\{\sigma_u,1/\sigma_l\}]$. Therefore, since $\sigma_{\min}(X)\sigma_{\min}(Y)\le\sigma_{\min}(XY)\le\sigma_{\max}(XY)\le\sigma_{\max}(X)\sigma_{\max}(Y)$,  it follows from \eqref{H0new} that the bounds for $N_k=S_kA_k^\#$ are 
\begin{equation}\label{Nbounds}
\lambda_l\min\{\sigma_l,1/\sigma_u\}\leq\sigma_{\min}(N_k)\le\sigma_{\max}(N_k)\le\lambda_u\max\{\sigma_u,1/\sigma_l\}.\end{equation}
Since $S_k$ is symmetric, $N_k^{\top}M_k=(A_k^{\#})^{\top}S_kS_k^{-1}=(A_k^{\#})^{\top}$ and $M_kN_k=S_k^{-1}S_kA_k^{\#}=A_k^{\#}$, so 
\begin{equation}\label{MNplus}\tfrac12(N_k^{\top}M_k+M_kN_k)=\tfrac12\big(A_k^{\#}+(A_k^{\#})^{\top}\big)=\tfrac12(A_k+A_k^\top) P+\tfrac12(A_k^{-1}+(A_k^{-1})^\top)P^\bot.\end{equation}
Let $0\neq x\in\H$ and $y=A_k^{-1}x$. From \eqref{H2new}, it holds that $\nu_l\|Px\|^2\le\tfrac12\langle(A_k+A_k^\top) P x,Px\rangle$, and since $P^\bot x=A_kP^\bot y$ implies $\|P^\bot x\|\le\sigma_u\|P^\bot y\|$, we obtain 
\begin{eqnarray*}
\frac{\nu_l}{\sigma_u^2}\|P^\bot x\|^2&\le&\nu_l\|P^\bot y\|^2\le\tfrac12\langle(A_k+A_k^\top) P^\bot  y,P^\bot y\rangle\\
&=&\langle A_k P^\bot  y,P^\bot y\rangle
=\langle P^\bot x,A_k^{-1}P^\bot x\rangle\\
&=&\tfrac12\langle(A_k^{-1}+(A_k^{-1})^\top) P^\bot  x,P^\bot x\rangle,
\end{eqnarray*}
so \eqref{MNplus} implies
\begin{equation}\label{NMbounds}
\nu_l\min\{1,1/\sigma_u^2\}\leq\lambda_{\min}(\tfrac12 (N_k^{\top}M_k+M_kN_k)).\end{equation}
Therefore, \eqref{Ncond} and \eqref{MNcond} follow from \eqref{Nbounds} and \eqref{NMbounds}, taking $\sigma_N^l:=\min\Big\{\lambda_l\min\{\sigma_l,\tfrac1{\sigma_u}\},\ \nu_l\min\{1,\tfrac1{\sigma_u^2}\}\Big\}$, and 
    $\sigma_N^u:=\lambda_u\max\{\sigma_u,\tfrac1{\sigma_l}\}$.
So statement (1) holds. 

Now we check that $(T^\#)^{-1}$ is calm at $(0,x^*+w^*)$. We consider in $\H^2=\H\times\H$ the norm given by $\|(x,w)\|=\sqrt{\|x\|^2+\|w\|^2}$ (the derivation is similar with any norm, since $\H^2$ is finite-dimensional). Since $\gph T$ and  $\SS\times\SS^\bot$ are subtransversal  at $(x^*,w^*)\in Z^*=\gph T\cap (\SS\times\SS^\bot)$, there exist constants $\delta>0$ and $L>0$ such that,  if $\|(x,w)-(x^*,w^*)\|\le\delta$,
\begin{equation}\label{subtransv1}
d((x,w),Z^*)\leq L\max \{d((x,w),\gph T),d((x,w), \SS\times\SS^\bot)\} .
\end{equation}
Since $Z^*\subset\SS\times\SS^\bot$, we have $d((x,w),Z^*)\ge d((x,w),\SS\times\SS^\bot)$, so we may assume without loss of generality that $L\geq1$.
Note that, for any $(x,w)\in\H^2$,  $P_{\SS\times\SS^\bot}(x,w)=(Px,P^\bot w)$,  where $P_{\SS\times\SS^\bot}$ is the orthogonal projection onto $\SS\times\SS^\bot$. Therefore, 
\begin{eqnarray}\label{distSS}\nonumber
d((x,w), \SS\times\SS^\bot)^2&=&\|(x,w)-(Px,P^\bot w)\|^2\\
        &=&\|x-Px\|^2+\|w-P^\bot w\|^2=\|P^\bot x\|^2+\|Pw\|^2.\label{Inecdist1}
\end{eqnarray}
Let $z\in\H$ and $v\in T^\#(z)$. Then, $Pv+P^\bot z\in T(Pz+P^\bot v)$. By defining $x=Pz+P^\bot v$ and $w=Pv+P^\bot z$, we obtain $(x,w)\in\gph T$, so 
\begin{equation}\label{distgT1}d((x,w),\gph T)=0.\end{equation}
Furthermore, $Px=Pz$, $P^\bot x=P^\bot v$, $Pw=Pv$ and $P^\bot w=P^\bot z$, so \eqref{distSS}
gives 
\begin{equation}\label{distStS1}d((x,w),\SS\times\SS^\bot)^2=\|P^\bot v\|^2+\|Pv\|^2=\|v\|^2.\end{equation}    
Similarly, let $(\bar x,\bar w)\in Z^*$, so $\bar z=\bar x+\bar w\in (T^\#)^{-1}(0)$. We obtain 
\begin{eqnarray}\nonumber
\|(x,w)-(\bar x,\bar w)\|^2&=&\|x-\bar x\|^2+\|w-\bar w\|^2\\&=&\|Px-\bar x\|^2+\|P^\bot x\|^2+\|Pw\|^2+\|P^\bot w-\bar w\|^2\nonumber \\
&=&\|Pz-P\bar z\|^2+\|P^\bot z-P^\bot\bar z\|^2+\|v\|^2=\|z-\bar z\|^2+\|v\|^2,\label{Equivnorm1}
\end{eqnarray}
so, taking the infimum in both sides,  
\begin{equation}\label{equivdist}
d((x,w),Z^*)^2=d(z,(T^\#)^{-1}(0))^2+\|v\|^2.\end{equation}
Let $U:=\{z\in\H: \|z-z^*\|\le \frac{\sqrt  2}{2}\delta\}$ and  $V:=\{v\in\H: \|v\|\le \frac{\sqrt  2}{2}\delta\}$. Let  $z\in U$ and suppose that $T^\#(z)\neq\emptyset$.  Then, for any $v\in T^\#(z)\cap V$, defining $x=Pz+P^\bot v$ and $w=Pv+P^\bot z$, it holds from \eqref{Equivnorm1} that $\|(x,w)-(x^*,w^*)\|=\sqrt{\|z-z^*\|^2+\|v\|^2}\le\delta,$
so \eqref{equivdist}, \eqref{subtransv1}, \eqref{distgT1} and \eqref{distStS1} give
\begin{eqnarray}\nonumber
    d(z,(T^\#)^{-1}(0))&=&\sqrt{d((x,w),Z^*)^2-\|v\|^2}\\\nonumber
                   &\le& \sqrt{L^2\max\{d((x,w),\gph T)^2,d((x,w),\SS\times\SS^\bot)^2\}-\|v\|^2}\\
                   &=&\sqrt{L^2-1}\,\|v\|.\label{subregdist1}
\end{eqnarray}
Since $v\in T^\#(z)\cap V$ is arbitrary, \eqref{subregdist1} implies 
\begin{equation}\label{subreg2}
 d(z,(T^\#)^{-1}(0))\le \sqrt{L^2-1}\;d(0, T^\#(z)\cap V).
\end{equation}
If  $T^\#(z)\cap V=\emptyset$, then $d(0, T^\#(z)\cap V)=\infty$ and \eqref{subreg2} holds trivially. 
Thus, $T^\#$ is metrically subregular at $z^*$ for 0, so $(T^\#)^{-1}$ is calm at 0 for $z^*=x^*+w^*$. 

With the identifications of Theorem \ref{thm:equiv}, the hypotheses in (a), (b) and (c) are exactly the same as in \cite[Theorem 4.4]{BLPS:25} (in particular, \eqref{H3new} is equivalent to \eqref{Mcond1}), so the results follow.

Finally, if $T$ is polyhedral so are $T^\#$ are  $(T^\#)^{-1}$. By the mentioned 
Robinson's theorem \cite{rob:81}, $(T^\#)^{-1}$ is calm at $0$ for every
$z^*\in Z^*$, so the results of item 2 apply.
\end{proof}

\begin{remark}
    The advantage of the polyhedral case is that subtransversality of $\gph T$ and $\SS\times\SS^\bot$ holds for every pair in $Z^*$, so we do not need to check it at the limit point $(x^*,w^*)$, which is a priori unknown. A classical example is $T=F+N_C$ with $F$ affine and $C$ convex polyhedral,which is the type of operators we will deal with in the numerical examples of section \ref{sec:numerics}.
\end{remark}

We now establish the relationship between Algorithm \ref{AIPDA} and the progressive decoupling scheme of \cite{rock:24} (and its specialization to the optimization case in \cite{Rock:26}.)
Basically, given $x_k\in \SS$ and $w_k\in\SS^\bot$, the iteration of the PDA in \cite[Section 4]{rock:24} chooses some self-adjoint positive definite linear operator $J_k$, solves
\begin{equation}\label{IPDA}
\x^{k}\approx (T^k)^{-1}(0), \text{ with } 
T^k(x)=T(x)+J_k(x-x_k)-w_k,
\end{equation}
and updates by $x_{k+1}=P\x^k$ and $w_{k+1}=w_k-J_kP^\bot \x^k.$
This procedure is equivalent to the variable metric proximal point method of \cite{rock:23b} aplied to the partial inverse of $T$, that is,
\begin{equation}\label{VMPPA}
x_{k+1}+w_{k+1}=:z_{k+1}\approx \tilde R_k(z_k),\, \text{ with }\, 
z_k:=x_k+w_k\, \text{ and }\,\tilde R_k:=(I+(J_k^\#)^{-1}T^\#)^{-1}.
\end{equation}
Assuming that \eqref{linkage} is solvable, if the sequence $\{J_k\}$ satisfies 
 \begin{eqnarray}\label{Jcond0}
&&0<\lambda_J^l\leq\lambda_{min}(J_k)\leq \lambda_{max}(J_k) \leq \lambda_J^u,\\ \label{Jcond1}
&&\frac{1}{1+\hat\zeta_k}J_k\preceq J_{k+1}\preceq(1+\hat\zeta_k)J_k,\quad \hat\zeta_k >0, \,
\quad \sum_{k=0}^{\infty} \hat\zeta_k<\infty,
\end{eqnarray}
and the sequence of errors is summable, 
the obtained sequence $\{(x_k,w_k)\}$ converges to some $(x^*,w^*)\in Z^*$. Also, under metric subregularity of $T^\#$ at $(x^*,w^*)$ for $0$, the distance $\dist_{J_k^\#}((x_k,w_k,Z^*)$ converges Q-linearly to zero if, for $e_k\ge 0$ with $\sum e_k<\infty,$\footnote{
\cite[Theorem 4.1]{rock:24} takes $\x^k$ verifying $2\dist_{J_k^{-1}}(0,T^k(\x^k))\le e_k\sqrt{\|x_{k+1}-x_k\|_{J_k}^2+\|w_{k+1}-w_k\|_{J_k^{-1}}^2}$, which implies \eqref{er:Rock}, a weaker condition used in \cite[Theorem 3.1]{rock:24}, which, in particular, does not require $\x^k\in\dom T^k$.}
\begin{equation}\label{er:Rock}
\|x_{k+1}+w_{k+1}-\tilde R_k(x_k+w_k)\|_{J_k^\#}\leq e_k\sqrt{\|x_{k+1}-x_k\|_{J_k}^2+\|w_{k+1}-w_k\|_{J_k^{-1}}^2}.
\end{equation}

It is proved in \cite[Theorem 5.1]{BLPS:25} that, for $k$ sufficiently large such that $e_k<1$, there exists a symmetric positive definite matrix $E_k$ that verifies
\begin{equation}\label{Econd}
    (J_k^\#)^{-\frac12}E_k(J_k^\#)^{\frac12}(z_k-\tilde R_k(z_k))=z_k-z_{k+1},
\end{equation}
with $E_k\to I$ as $k\to\infty$.  Then, by defining
\begin{equation}\label{MNprox}
    M_k=(J_k^\#)^{\frac12}E_k^{-1}(J_k^\#)^{\frac12}\quad\text{and}\quad
     N_k=(J_k^\#)^{-\frac12}E_k(J_k^\#)^{\frac12},
\end{equation}
and using the parameters $\sigma_k\equiv 0$, $\varepsilon_k\equiv 0$, $c_k\equiv 1$, it can be proved that, from some iteration onwards, the algorithm \eqref{VMPPA} with tolerance condition \eqref{er:Rock} can be cast into the framework of the scheme \eqref{hippm}-\eqref{Update1}. The issue here is that condition \eqref{er:Rock} is not sufficient to guarantee that $E_k$ preserves the subspaces $\SS$ and $\SS^\bot$, which is necessary for specialization to algorithm \ref{AIPDA}, and therefore we must formulate a stronger assumption. We require condition \eqref{er:Rock} to be satisfied in each block, that is, 
\begin{equation}\label{er:Rock1}\begin{array}{rcl}
\|x_{k+1}-P\tilde R_k(x_k+w_k)\|_{J_k}&\leq& e_k\|x_{k+1}-x_k\|_{J_k},\\
\|w_{k+1}-P^\bot\tilde R_k(x_k+w_k)\|_{J_k^{-1}}&\leq& e_k\|w_{k+1}-w_k\|_{J_k^{-1}}.\end{array}
\end{equation}
The next theorem shows that, for $k$ large enough, any pair $x_{k+1}\in\SS$, $w_{k+1}\in\SS^\bot$  satisfying \eqref{er:Rock1} can be obtained by a particular iteration of Algorithm \ref{AIPDA}, so it inherits some of its convergence properties.

\begin{theorem}\label{rel-sum}
Let $\{(x_k,w_k)\}$ be a sequence generated by \eqref{IPDA}, with self-adjoint positive-definite operators $J_k$ preserving $\SS$ and $\SS^\bot$, satisfying \eqref{er:Rock1} with $e_k\ge0$ and $\sum_k e_k<\infty$. Then there exists $K\in\mathbb{N}$ such that, for $k\ge K$, the iterates of Algorithm~\ref{AIPDA} started at $(x_K,w_K)$ with $A_k=J_k$, $c_k\equiv1$, $\sigma_k\equiv0$, $\varepsilon_k\equiv0$ and suitable $S_k$ and $\tau_k\in[1-\theta,1+\theta]$ coincide with $\{(x_k,w_k)\}$. If, in addition, $\{J_k\}$ satisfies \eqref{Jcond0}-\eqref{Jcond1}, then $\{(x_k,w_k)\}$ converges to some $(x^*,w^*)\in Z^*$ and, if $\gph T$ and $\SS\times\SS^\bot$ are subtransversal at $(x^*,w^*)$, the convergence is R-linear in the norm induced by $(J^\#_\infty)^{-1}$.
\end{theorem}
 \begin{proof}
    The proof follows the ideas of \cite[Theorem 5.1 and Corollary 5.2]{BLPS:25}, so we give the general sketch. Assume that $k$ is large enough so that $0\le e_k<1$. 
If $e_k =0$ then the subproblem is solved exactly, wich can be done by Algorithm \ref{AIPDA} with $A_k:=J_k$, $S_k:=(J_k^\#)^{-1}$ and $\tau_k=1$.

	 Let $e_k >0$, $(x_k,w_k)$ the current iterate in
	 \eqref{IPDA}, $(x_{k+1},w_{k+1})$ any update satisfying \eqref{er:Rock1}, $x^*_k=P\tilde R_k(x_k+w_k)$ and $w^*_k=P^\bot\tilde R_k(x_k+w_k)$. 
Let $J_{k,\SS}$ and $J_{k,\SS^\bot}$ be the blocks of $J_k$ corresponding to $\SS$ and $\SS^\bot$, respectively, so $J_k^\#=J_{k,\SS}P+J_{k,\SS^\bot}^{-1}P^\bot$. 

If $x_{k+1}=x_k$, then \eqref{er:Rock1} gives $x^*_k=x_k$ and we take $E_{k,\SS}=I_\SS$ (analogously in $\SS^\bot$ if $w_{k+1}=w_k$); hence we may assume that $x_{k+1}\neq x_k$.
For $d_k:=\|J_{k,\SS}^{\frac12}(x_k - x_{k+1})\|$, there exists an orthogonal matrix $C_k:\SS\to\SS$ such that 
$C_kJ_{k,\SS}^{\frac12}(x_k-x_{k+1})=(d_k,0,...,0)$ and $C_kJ_{k,\SS}^{\frac12}(x_k-x^*_k)=(y_1^k,y_2^k,0,...,0),$ so
$$C_kJ_{k,\SS}^{\frac12}(x_{k+1}-x^*_k)=(y_1^k-d_k,y_2^k,0,...,0)$$ and
$$ \|(y_1^k-d_k,y_2^k,0,...,0)\|=\|x_{k+1} -x^*_k\|_{J_{k,\SS}}\\
\leq e_k\|x_k -x_{k+1}\|_{J_{k,\SS}}=e_k d_k.$$
Hence,  $|y_1^k-d_k|\leq e_kd_k$ and $|y_2^k|\leq e_kd_k$, which implies
$0<(1-e_k)d_k\leq y_1^k\leq (1+e_k)d_k.$ (Note that this can not be achieved if the iterates satisfy \eqref{er:Rock} but not \eqref{er:Rock1}.)
Defining $\hat a_k= d_k/y_1^k + (y_2^k)^2/(y_1^k)^2$ and $\hat b_k= -y_2^k/y_1^k$,  we observe that  $\hat a_k>0$ and $\hat a_k-\hat b_k^2=d_k/y_1^k>0$, so a straightforward computation shows that the matrix  
$$\bar E_{k,\SS}:=\left(\begin{array}{cccc}
    \hat a_k & \hat b_k &  0    &  0 \\
    \hat b_k &    1     &  0    &  0 \\
        0    &    0     & \ddots&  0\\ 
         0    &    0     &    0 &  1\\
\end{array}\right)$$
is positive definite and verifies $\bar E_{k,\SS}(y_1^k,y_2^k,0,...,0)=(d_k,0,...,0)$. Also,
\begin{equation}\label{cotasmat}
0\leq\frac{|y_2^k|}{y_1^k}\leq\frac{e_k}{1-e_k}\quad\text{ and }\quad
\frac{1}{1+e_k}\leq \frac{d_k}{y_1^k}\leq\frac{1}{1-e_k},
\end{equation}
so $\|\bar E_{k,\SS}-I_\SS\|=O(e_k)$. By defining  the self-adjoint positive-definite linear operator $E_{k,\SS}:\SS\to\SS$ as $E_{k,\SS}:=C_k^{-1}\bar E_{k,\SS}C_k$, we obtain that$J_{k,\SS}^{-\frac12} E_{k,\SS} J_{k,\SS}^{\frac12} (x_k-x^*_k)= x_k-x_{k+1}$.
     
Analogously, there exists a self-adjoint positive-definite linear operator $E_{k,S^\bot}:\SS^\bot\to\SS^\bot$ such that $\|E_{k,\SS^\bot}- I_{\SS^\bot}\|=O(e_k)$ and
 $J_{k,\SS^\bot}^{\frac12} E_{k,\SS^\bot} J_{k,\SS^\bot}^{-\frac12} (w_k-w^*_k)= w_k-w_{k+1}$, so we obtain that the self-adjoint positive-definite linear operator $E_k$, given by $E_k(x,w):=(E_{k,\SS}\ x,{E_{k,\SS^\bot}\ w})$, preserves the subspaces $\SS$ and $\SS^\bot$, satisfies \eqref{Econd}, and verifies $\|E_k-I\|=O(e_k)$. Hence, \begin{equation}\label{Eksum}
     \sum_k\|E_k-I\|<\infty,
 \end{equation} 
and, in particular, $E_k\to I$.
Note that $x^*_k+w^*_k=\tilde R_k(x_k+w_k)$ is equivalent to
$$(J_k^\#)^{\frac12}(x_k-x^*_k+w_k-w^*_k) \in (J_k^\#)^{-\frac12}T^\#(x^*_k+w^*_k),$$
so multiplying this inclusion by  $ (J_k^\#)^{-\frac12}E_k$, choosing $c_k=1$, $M_k, N_k$ as in \eqref{MNprox} and  defining $\hat v^k=N_k(x_k-x^*_k+w_k-w^*_k)$, we obtain the system \eqref{hippm} with $\delta_k=\varepsilon_k=0$, that is,
\begin{equation}\label{Spezhippm}
    \left\{
\begin{array}{lcl}
\hat v^k &\in& M_k^{-1}T^\#(x^k),\\ 
0&=& \hat v^k+N_k(x^*_k-x_k+w^*_k-w_k)=\hat v^k+x_{k+1}-x_k+w_{k+1}-w_k,
\end{array}
\right.\end{equation}
which trivially satisfies the tolerance 
 conditions in \eqref{hippm} with $\sigma_k=0$. Note that both $M_k$ and 
$N_k$ are positive definite and preserve the subspaces $\SS$ and $\SS^\bot$, with $M_k$ symmetric and $N_k$ not necessarily symmetric. Also, $N_k^TM_k=J^\#_k$ is positive definite and, since $E_k\to I_n$, we have that $M_k\to J^\#_\infty$ and $N_k\to I_n$. 
A simple computation  in \eqref{Update1} gives 
$$a_k=\frac {\|
(J^\#_k)^{\frac12}( x_k-x^*_k+w_k-w^*_k)\|^2} {\|
(J^\#_k)^{\frac12}( x_k-x^*_k+w_k-w^*_k)\|^2_{E_k}},
$$
and, since $E_k\to I$, it holds that 
$1-\theta\leq \lambda_{min}(E_k)\leq\lambda_{max}(E_k)\leq 1+\theta$ from some  $ K$ onwards. 
Therefore by choosing 
$\tau_k=a_k^{-1}$, we have $\tau_k\in[1-\theta,1+\theta]$ and, using \eqref{Spezhippm}, we obtain  
$\bar z^{k+1}= z^k-\tau_k a_k \hat v^k = x_k+w_k-\hat v^k = x_{k+1}+w_{k+1}.$

In view of Theorem \ref{thm:equiv}, the scheme \eqref{IPDA} with tolerance condition \eqref{er:Rock1} can be cast in the framework of Algorithm \ref{AIPDA} if the operators $A_k:=J_k$ and  $S_k:=M_k^{-1}$ (note that both are symmetric positive-definite and preserve subspaces $\SS$ and $\SS^\bot$), satisfy $N_k=S_kA_k^\#$. From \eqref{MNprox}, we have
$$S_kA_k^\#=(J_k^\#)^{-\frac12}E_k(J_k^\#)^{-\frac12}J_k^\#=(J_k^\#)^{-\frac12}E_k(J_k^\#)^{\frac12}=N_k,$$
so the result follows. 

Finally, under \eqref{Jcond0}-\eqref{Jcond1}, from the definition of $M_k$ and \eqref{Eksum}, we obtain that  $S_k=M_k^{-1}$ satisfies \eqref{H3new}. Moreover, $J_k\to J_\infty$ and $S_k\to(J^\#_\infty)^{-1}$. Therefore statement 2(c) of Corollary~\ref{ConvProp} applies under the subtransversality assumption, and the claimed convergence holds.\end{proof}

\section{AIPHA: progressive hedging for stochastic variational inequalities}\label{sec:aipha}

We now specialize AIPDA to multistage stochastic variational inequalities, in the formulation given by Rockafellar and Wets (see \cite{RW:17} for a detailed exposition). Consider $N$ stages: at each stage $t=1,\dots,N$, a decision $x_t\in\R^{n_t}$ is taken and then an observation $\xi_t\in \Xi_t$ is revealed, where each $\Xi_t$ is finite. Hence the set $\Xi$ of scenarios $\xi=(\xi_1,\dots,\xi_N)$ is finite and, assuming that each scenario has a known probability $p(\xi)>0$ with $\sum_{\xi\in\Xi}p(\xi)=1$, $\Xi$ is a finite probability space. Interpreting decisions as responses to the information, we consider the linear space
\begin{equation}\label{fHs}
   \L=\{x :\Xi\to\R^n=\R^{n_1}\times\dots\times\R^{n_N},\ \xi\mapsto (x_1(\xi),\dots,x_N(\xi))\}.
\end{equation}
If $\Xi$ has $s$ elements, $\L$ can be identified with $s$ copies of $\R^n$, a linear space of dimension $ns$. Equipped with the expectational inner product
$$\langle x,w\rangle=\sum\limits_{\xi\in\Xi}p(\xi) \left(x(\xi), w(\xi)\right),$$
where $(\cdot,\cdot)$ stands for the usual inner product in $\R^n$, it becomes a Hilbert space, whose norm is denoted by $\|\cdot\|$. Since the decision at stage $t$ can depend only on the information revealed before it, that is, $x_t(\xi)=x_t(\xi_1,\dots,\xi_{t-1})$, the concept of \emph{nonanticipativity} comes into play, and $\L$ decomposes into the orthogonal subspaces $\NN$ and $\M$, respectively the \emph{nonanticipativity subspace} and the \emph{multipliers subspace}, given by
\begin{equation}\label{NAsubs}
   \NN=\{x \in\L: x_t(\xi) \text{ does not depend on } \xi_t,\dots,\xi_N,\ t=1,\dots,N\}
\end{equation}
and
\begin{equation}\label{Msubs}
   \M=\{w \in\L:\E[w_t(\xi)\mid\xi_1,\dots,\xi_{t-1}]=0,\ t=1,\dots,N\}.
\end{equation}
Having functions $F(\cdot,\xi):\R^n\to\R^n$ for each $\xi\in\Xi$, the mapping $\F:\L\to\L$ is obtained by defining $\F(x )$ as the application $\Xi\ni\xi\mapsto F(x(\xi),\xi)\in\R^n$, which is continuous and monotone if the functions $F(\cdot,\xi)$ are so. Considering also the constraint set
\begin{equation}\label{Cset}
   \C=\{x \in\L:x(\xi)\in C(\xi),\ \forall \xi\in\Xi\},
\end{equation}
where the sets $C(\xi)$ are closed and convex, $\C$ is also closed and convex. Then the operators $\F+N_\C$ and $\F+N_{\N\cap \C}$ (the latter when $\N\cap\C\neq\emptyset$) are maximally monotone.
The \emph{stochastic variational inequality in basic form} is stated as the problem of finding $x \in\C$ such that
\begin{equation}\label{SVI1}
   0\in \left(\F+N_{\N\cap \C}\right)(x ).
\end{equation}
Under appropriate constraint qualifications, the inclusion \eqref{SVI1} is equivalent to the \emph{stochastic variational inequality in extensive form}\footnote{In \cite{RS:18} the extensive form is written as $-\F(x)-w\in N_\C(x)$. Both forms are equivalent through $w\mapsto -w$, since $\M$ is a subspace; we use \eqref{SVI2} throughout, also in Section~\ref{sec:numerics}.}
\begin{equation}\label{SVI2}
   w \in \left(\F+N_\C\right)(x ),\quad x \in\NN,\;w
\in \M,
\end{equation}
which is a linkage problem of the type \eqref{linkage} with $\H=\L$, $\SS=\NN$ and $\SS^\bot=\M$, so Algorithm \ref{AIPDA} applies to $T=\F+N_\C$.

The resulting method is of progressive hedging type: as we now explain, its proximal step decomposes by scenarios. For each $\xi\in\Xi$, let $\L^\xi$ be the copy of $\R^n$ formed by the $\xi$-coordinates of $\L$, and let $T_\xi=F(\cdot,\xi)+N_{C(\xi)}$. The operator $T=\F+N_\C$ is \emph{separable}, in the sense that $v\in T(x)$ if and only if $v(\xi)\in T_\xi(x(\xi))$ for all $\xi\in\Xi$. Moreover, if $\varepsilon(\xi)\ge0$ for all $\xi$ and $\varepsilon=\sum_{\xi\in\Xi}p(\xi)\varepsilon(\xi)$, then for any $(w,u)\in\gph T$,
\[\begin{array}{rcl}
    \left(v(\xi)-u(\xi),z(\xi)-w(\xi)\right) &\geq& -\varepsilon(\xi),\;\forall \xi\in\Xi \quad \Rightarrow \\
   \Rightarrow   \langle v-u,z-w\rangle&=&\sum\limits_{\xi\in\Xi}p(\xi)\left(v(\xi)-u(\xi),z(\xi)-w(\xi)\right) \geq -\varepsilon,
\end{array}
\]
so that, using also \eqref{enlargesum},
\begin{equation}\label{EnSce}
    v(\xi)\in F(\cdot,\xi)^{\varepsilon(\xi)}(z(\xi))+N_{C(\xi)}(z(\xi))\subset T_\xi^{\varepsilon(\xi)}(z(\xi)),\; \forall \xi\in\Xi\ \Longrightarrow\ v\in T^{\varepsilon}(z).
\end{equation}
In contrast, the projection $P=P_\NN$ is \emph{not} separable by scenarios, since it averages the decisions of each stage over the scenarios that share the same history. Consequently, the operator $D_k=c_k^{-1}P+c_kP^\bot$ in \eqref{ips} couples the scenarios unless $c_k=1$, in which case $D_k=I$. Therefore, we take $c_k\equiv1$ and we require $A_k$ to act scenario-wise, that is, $(A_kx)(\xi)=A_k^\xi x(\xi)$, besides preserving $\NN$ and $\M$. Then, by \eqref{EnSce}, the inclusion in \eqref{ips} is implied by the $s$ independent inclusions
\begin{equation}\label{ips-ph}
A_k^\xi\big(x_k(\xi)-\x^k(\xi)\big)+w_k(\xi) \in (F(\cdot,\xi)^{\varepsilon_k(\xi)}+ N_{C(\xi)})(\w^k(\xi)),\qquad \xi\in\Xi,
\end{equation}
with $\varepsilon_k=\sum_{\xi}p(\xi)\varepsilon_k(\xi)$. The remaining quantities in \eqref{eq:conddeltaeps1}, \eqref{eq:conddeltaeps2} and \eqref{updates} involve $P\x^k$, $P\w^k$ and norms induced by $S_k$; they are global, but they are computed from the scenario solutions through a single averaging step (conditional expectations), so $S_k$ need not act scenario-wise (in Section~\ref{sec:numerics}, $S_k=P+\kappa P^\bot$).

\begin{remark}\label{rem:twostage}
For two-stage problems ($N=2$, $x=(x_1,x_2)$) with $s\ge2$, a scenario-wise operator $A_k$ preserves $\NN$ and $\M$ if and only if $A_k^\xi=\mathrm{blkdiag}\big(A_{k,1},A_{k,2}(\xi)\big)$ with $A_{k,1}$ independent of $\xi$. In particular, the choice $A_k^\xi=\Lambda-M_\xi$ used in Section~\ref{sec:numerics}, which gives, for each scenario, diagonal subproblems, requires the Jacobian $M_\xi$ of $F(\cdot,\xi)$ to have a first-stage block that is independent of the scenarios and no inter-stage coupling. Both properties hold in all the models of Section~\ref{sec:numerics}.
\end{remark}

With the preceding considerations, AIPDA leads to the following algorithm.

 \begin{algorithm}[H]
 \caption{Asymmetric Inexact Progressive Hedging Algorithm (AIPHA)}\label{AIPHA-eps}
\begin{algorithmic}[1]

\item \textit{Initialization:} Choose $x_0\in \NN$, $w_0\in \M$, $\bar{\sigma} \in (0,1)$ and $\theta \in (0,1)$.

\item \textit{Scenario proximal step:} Choose $\sigma_k\in [0,\bar{\sigma})$, a self-adjoint positive-definite linear mapping $S_k$ and a positive-definite linear mapping $A_k$, both preserving $\NN$ and $\M$, with $A_k$ acting scenario-wise, $(A_kx)(\xi)=A_k^\xi x(\xi)$. For each $\xi \in \Xi$, find $\x^k(\xi)$, $\w^k(\xi)$ and $\varepsilon_k(\xi)\geq 0$ satisfying \eqref{ips-ph}.

\item \textit{Tolerance test:} Set $\varepsilon_k=\sum_{\xi\in\Xi}p(\xi)\varepsilon_k(\xi)$, $\hat\delta^k=\w^k-\x^k$, compute $P\x^k$, $P\w^k$ and check
\begin{equation} \begin{array}{l}
\|(A_kP+P^\bot)\hat \delta^k\|_{S_k}^2 +2\varepsilon_k  \leq  \\\qquad\; \sigma_k^2\big(\| A_k(x_k-P\x^k)+P^\bot \w^k\|_{S_k}^2+\| A_k(x_k-P\w^k)+P^\bot \x^k\|_{S_k}^2\big),
\label{eq:conddeltaeps3}\end{array}
\end{equation}
\begin{equation} \begin{array}{l}
\langle \hat\delta^k,A_k^\top(P\w^k-x_k)-A_k P^\bot\x^k\rangle+\varepsilon_k\leq \\\qquad\;
\sigma_k^2\langle A_k(P\w^k-x_k+P^\bot\x^k),P\w^k-x_k+P^\bot\x^k\rangle.
\label{eq:conddeltaeps4}\end{array}
\end{equation}
If \eqref{eq:conddeltaeps3} or \eqref{eq:conddeltaeps4} fails, solve the scenario subproblems more accurately and repeat this step.

\item \textit{Update:} If $x_k-P\w^k=P^\bot \x^k=0$, STOP. Otherwise, choose  $\tau_k\in [1-\theta,1+\theta]$ and set
 \begin{equation}\label{updates1}
 \begin{array}{lcl}
 x_{k+1}&=& x_k - \tau_k \alpha_k S_kA_k(x_k-P\x^k)
\\
 w_{k+1}&=&w_k-\tau_k \alpha_k S_kP^\bot\w^k,
 \end{array}
 \end{equation}
where
\begin{equation*}
    \alpha_k
   = \frac{\langle A_k(x_k-P\x^k+P^\bot\x^k),x_k -P\w^k + P^\bot \w^k\rangle - \varepsilon_k} {\|A_k(x_k-P\x^k)+P^\bot \w^k\|_{S_k}^2}.
\end{equation*}
Set $k:=k+1$ and go to the Scenario proximal step.
\end{algorithmic}
\end{algorithm}

Algorithm~\ref{AIPHA-eps} is Algorithm~\ref{AIPDA} with $c_k\equiv1$, so Corollary \ref{ConvProp} applies (since $c_k\equiv1$, only the linear-rate part of statement 2(b) is available). In particular, if each $F(\cdot,\xi)$ is affine and monotone and each $C(\xi)$ is convex polyhedral, as in Section~\ref{sec:numerics}, AIPHA converges globally and, under the parameter choices of Corollary~\ref{ConvProp}, at a local linear rate.

With the particular choices $A_k\equiv rI$ and $S_k\equiv r^{-2}P+P^\bot$, for some $r>0$, and $\varepsilon_k(\xi)\equiv 0$, for all $\xi$, the inclusions \eqref{ips-ph} turn into
\begin{equation}\label{ips-ph-viejo}
r\big(x_k(\xi)-\x^k(\xi)\big)+w_k(\xi) \in (F(\cdot,\xi)+ N_{C(\xi)})(\w^k(\xi)),\qquad \xi\in\Xi,
\end{equation}
and inequality \eqref{eq:conddeltaeps3} becomes 
\begin{equation}
\|\hat\delta^k\|^2  \leq \sigma_k^2\big(\|x_k-P\x^k+P^\bot \w^k\|^2+\| x_k-P\w^k+P^\bot \x^k\|^2\big).
\label{eq:conddelta-viejo}\end{equation}
Since $\hat\delta^k=\w^k-\x^k=(x_k-P\x^k+P^\bot\w^k)-(x_k-P\w^k+P^\bot\x^k)$, the above inequality gives
$$\begin{array}{l}\frac{1-\sigma_k^2}{2}\left(\|x_k-P\x^k+P^\bot\w^k\|^2+\|x_k-P\w^k+P^\bot\x^k\|^2\right)\le\\\qquad\; \langle x_k-P\x^k+P^\bot\w^k,x_k-P\w^k+P^\bot\x^k\rangle\\\qquad\;\le \|x_k-P\x^k+P^\bot\w^k\|\|x_k-P\w^k+P^\bot\x^k\|,\end{array}$$
which, by  denoting $t:=\|x_k-P\w^k+P^\bot\x^k\|$, is equivalent the quadratic inequality
$$\frac{1-\sigma_k^2}{2}t^2-t\|x_k-P\x^k+P^\bot\w^k\|+\frac{1-\sigma_k^2}{2}\|x_k-P\x^k+P^\bot\w^k\|^2\le 0, $$
whose solutions give
\begin{equation}\label{des:ph-viejo}
    \frac{1-\rho_k}{1-\sigma_k^2}\|x_k-P\x^k+P^\bot\w^k\|\le\|x_k-P\w^k+P^\bot\x^k\|\le\frac{1+\rho_k}{1-\sigma_k^2}\|x_k-P\x^k+P^\bot\w^k\|,
\end{equation}
for $\rho_k:=\sqrt{1-(1-\sigma_k^2)^2}$. Therefore, using \eqref{eq:conddelta-viejo}, \eqref{des:ph-viejo} and $\rho_k\in[0,1)$, it holds that
\begin{align*}
    &\langle \hat\delta^k,P\w^k-x_k-P^\bot\x^k\rangle
    \\&\quad=\frac12\Big(\|\hat\delta^k\|^2+\|P\w^k-x_k-P^\bot\x^k\|^2
    \\&\qquad\quad-\|\hat\delta^k-P\w^k+x_k+P^\bot\x^k\|^2\Big)\\
    &\quad=\frac12\left(\|\hat\delta^k\|^2+\|P\w^k-x_k-P^\bot\x^k\|^2-\|x_k-P\x^k+P^\bot\w^k\|^2\right)\\
    &\quad\le\frac{\sigma_k^2+1}{2}\|x_k-P\w^k+P^\bot\x^k\|^2+\frac{\sigma_k^2-1}{2}\|x_k-P\x^k+P^\bot\w^k\|^2\\
    &\quad\le\left(\frac{(1-\sigma_k)^2(\sigma_k^2-1)}{2(1-\rho_k)^2}+\frac{\sigma_k^2+1}{2}\right)\|x_k-P\w^k+P^\bot\x^k\|^2\\
    &\quad\le\sigma_k^2\|x_k-P\w^k+P^\bot\x^k\|^2,
\end{align*}
which, with the chosen parameters, coincides with \eqref{eq:conddeltaeps4}, that turns to be a superfluous requirement.
Concerning the updates, it holds that $\alpha_k=r\hat\alpha_k$, where
$\hat\alpha_k=\frac{\langle x_k-P\x^k+P^\bot\w^k,x_k-P\w^k+P^\bot\x^k\rangle}{\|\langle x_k-P\x^k+P^\bot\w^k\|^2},$
so the updates are
 \begin{equation}\label{updates2}
 \begin{array}{lcl}
 x_{k+1}&=& x_k - \tau_k \hat\alpha_k (x_k-P\x^k),
\\
 w_{k+1}&=&w_k-\tau_k r\hat\alpha_k P^\bot\w^k.
 \end{array}
 \end{equation}
In summary, the algorithm solves \eqref{ips-ph-viejo} until condition \eqref{eq:conddelta-viejo} is achieved and updates using \eqref{updates2}. Through the identification $w_k\mapsto -w_k$, this is the inexact progressive hedging algorithm of \cite{BLP:24}, which in turn extends the exact progressive hedging algorithm of \cite{RS:18}, so our scheme encompasses both methods.

\section{Numerical experiments: electricity markets}\label{sec:numerics}

All experiments concern a family of two-stage stochastic Cournot games arising in
electricity markets, organized so that the \emph{asymmetry} of the operator increases
from zero in a controlled and economically meaningful way; this isolates the role of
the (possibly nonsymmetric) metric $A$ and of the outer update of
Algorithm~\ref{AIPHA-eps}. Throughout, $F(z,\xi)=M_\xi z+b_\xi$, $z(\xi)\in C(\xi)$,
$\NN$ is the nonanticipativity subspace, $\MM=\NN^\perp$, and we measure asymmetry by
\begin{equation}\label{eq:eta}
\eta \;=\; \max_{\xi}\;
\frac{\norm{\skw(M_\xi)}_2}{\norm{\sym(M_\xi)}_2},
\qquad \skw(M)=\tfrac12(M-M^\top),\quad \sym(M)=\tfrac12(M+M^\top).
\end{equation}

\subsection{A nested family of market models}\label{sec:family}

\paragraph{P1: symmetric benchmark.}
Four producers choose first-stage quantities $x_i$ and recourse quantities
$y_i(\xi)$, with inverse demands $p^1(x)=\alpha^1(a^1-\sum_i x_i)$ and
$p^2(y,\xi)=\alpha^2_\xi(a^2_\xi-\sum_i y_i)$, costs
$C_i=c^1_ix_i-\gamma_i p^1x_i+c^2_i(\xi)y_i-\gamma_i p^2y_i$, and constraints
$x_i,y_i\ge0$, $x_i+y_i\le\ell_i$. With equal discount factors ($\gamma_i\equiv\gamma$)
the Jacobian is symmetric ($\eta=0$) and the problem is a two-stage stochastic program.

\paragraph{P2: nonsymmetric market power.}
Player-specific discount factors $\gamma=(1,1,0.7,0.7)$ and a \emph{directional}
cross-market coupling $\delta\ge0$ (group~2 affects the price of group~1 with weight
$1+\delta$, and conversely with weight $1-\delta$) give the per-stage matrix
$M_{\mathrm{st}}=\alpha\bigl(\begin{smallmatrix}\gamma_1(J+I)&(1+\delta)\gamma_1J\\
(1-\delta)\gamma_2J&\gamma_2(J+I)\end{smallmatrix}\bigr)$, $J=\mathbf{1}\mathbf{1}^\top$.
Sweeping $\delta\in[0,8]$ yields $\eta\in[0.07,2.08]$. The game is monotone only up to
$\delta\approx2.7$ ($\eta\approx0.96$); $\delta\in\{4,8\}$ is kept as a stress test,
outside the scope of the theory.

\paragraph{P3: priced network congestion.}
The producers share a transmission corridor whose capacity constraints are
\emph{priced} rather than inserted into $C(\xi)$, which would destroy its
separability (Section~\ref{sec:P3detail}). The extended Jacobian is
\begin{equation}\label{eq:Mred}
M_\xi=\begin{pmatrix} G_\xi & E\\ -E^\top & 0\end{pmatrix},
\qquad E=\mathrm{diag}(e,e),
\end{equation}
a \emph{purely antisymmetric} price coupling: $\sym(M_\xi)$ and monotonicity are
unaffected, while $\eta$ grows with the congestion scale $s$ of the injection factors
$e$. The first-stage network constraint is active at the solution for every $s$, and
the second-stage one in two of the five scenarios.

\paragraph{P4: nodal markets of increasing size.}
$n$ producers form local Cournot markets of four, coupled \emph{only} through the
corridor. Then $\norm{\sym(G)}=O(1)$ while $\norm{E}\sim\sqrt{n}$, so $\eta\sim\sqrt{n}$:
large realistic systems live in the high-asymmetry regime.

\subsection{The congested market model in detail}\label{sec:P3detail}

Producers $i=1,\dots,n$ are partitioned into local markets $\ell=1,\dots,L$ of size
$m$ ($m=n=4$ in P3). The random data $\xi\mapsto(\alpha^2_\xi,a^2_\xi,c^2(\xi),T^2(\xi))$
live on a finite set $\Xi$ with probabilities $p_\xi$; first-stage data
$(\alpha^1,a^1,c^1,T^1)$ are deterministic. Local inverse demands are
$p^1_\ell(x)=\alpha^1(a^1-\sum_{j\in\ell}x_j)$ and
$p^2_\ell(y,\xi)=\alpha^2_\xi(a^2_\xi-\sum_{j\in\ell}y_j)$. Given its rivals' decisions,
producer $i$ (in market $\ell(i)$, with discount factor $\gamma_i$) solves
\begin{equation}\label{eq:player}
\begin{aligned}
\min_{x_i,\;y_i(\cdot)}\quad
& c^1_ix_i-\gamma_i\,p^1_{\ell(i)}(x)\,x_i
+\sum_{\xi\in\Xi}p_\xi\Bigl[c^2_i(\xi)y_i(\xi)
  -\gamma_i\,p^2_{\ell(i)}(y(\xi),\xi)\,y_i(\xi)\Bigr]\\
\text{s.t.}\quad
& x_i\ge0,\quad y_i(\xi)\ge0,\quad x_i+y_i(\xi)\le\ell_i \quad(\xi\in\Xi),\\
& e^\top x\le T^1,\qquad e^\top y(\xi)\le T^2(\xi)\quad(\xi\in\Xi),
\end{aligned}
\end{equation}
where $e\in\mathbb{R}^n$ collects the signed injection factors of the corridor. The
last line couples all players, so \eqref{eq:player} is a generalized Nash game. We
compute its \emph{variational equilibria}, characterized by shared-constraint
multipliers common to all players \cite{fk:10}: congestion prices $\nu\ge0$ and
$\mu(\xi)\ge0$, under which producer $i$ pays $e_i\nu\,x_i+e_i\mu(\xi)\,y_i(\xi)$.
Stationarity and complementarity give, for $z(\xi)=(x,y(\xi),\nu,\mu(\xi))\in\mathbb{R}^{2n+2}$,
\begin{equation}\label{eq:Fdetail}
F(z,\xi)=
\begin{pmatrix}
\Bigl[c^1_i+\gamma_i\alpha^1\bigl(\sum_{j\in\ell(i)}x_j+x_i-a^1\bigr)
      +e_i\nu\Bigr]_{i=1}^{n}\\[4pt]
\Bigl[c^2_i(\xi)+\gamma_i\alpha^2_\xi\bigl(\sum_{j\in\ell(i)}y_j+y_i-a^2_\xi\bigr)
      +e_i\mu\Bigr]_{i=1}^{n}\\[4pt]
T^1-e^\top x\\[2pt]
T^2(\xi)-e^\top y
\end{pmatrix},
\qquad
C(\xi)=\prod_{i=1}^{n}\Delta_i\;\times\;\mathbb{R}_+^{2},
\end{equation}
with $\Delta_i=\{(x_i,y_i):x_i,y_i\ge0,\ x_i+y_i\le\ell_i\}$, so $C(\xi)$ is separable.
With $\NN=\{z(\cdot):x(\xi)\equiv\bar x,\ \nu(\xi)\equiv\bar\nu\}$ (the first-stage
price is also here-and-now) and $\MM=\NN^\perp$ under
$\langle u,v\rangle=\sum_\xi p_\xi\,u(\xi)^\top v(\xi)$, the variational equilibria
are exactly the solutions of the stochastic variational inequality: find $z\in\NN$,
$w\in\MM$ with
\begin{equation}\label{eq:SVI}
w(\xi)-F(z(\xi),\xi)\in N_{C(\xi)}(z(\xi)),\qquad \xi\in\Xi .
\end{equation}
Here $F(\cdot,\xi)$ is affine with $M_\xi$ as in \eqref{eq:Mred},
$G_\xi=\mathrm{blkdiag}_\ell\bigl(\Gamma_\ell\alpha(J_m+I_m)\bigr)$ per stage, and
$\langle M_\xi z,z\rangle=\langle \sym(G_\xi)(x,y),(x,y)\rangle\ge0$, since what
producers pay is what the constraint collects (we verify
$\lambda_{\min}(\sym G_\xi)\ge0$ on every instance). Each
$K(\xi)=\{(x,y)\in\prod_i\Delta_i:\ e^\top x\le T^1,\ e^\top y\le T^2(\xi)\}$ is a
nonempty compact polyhedron, so variational equilibria and congestion prices exist
without Slater-type conditions. $F(\cdot,\xi)$ is not strongly monotone (the price
block of $\sym(M_\xi)$ vanishes), but the subproblems of all methods are
($M_\xi+rI$ in progressive hedging, $M_\xi+A=\Lambda\succ0$ in AIPHA), so no
regularization is needed.

\paragraph{Data.}
$n=4$ (P3) or $n\in\{4,\dots,192\}$ with $m=4$ (P4); $|\Xi|=5$ equiprobable
scenarios; $\gamma_i\sim U(0.65,1)$; $c^1_i\sim U(0.5,1.5)$, $c^2_i(\xi)\sim U(0.5,2)$;
$\alpha^1=1$, $\alpha^2_\xi\sim U(0.8,1.6)$; $a^1=10$, $a^2_\xi\sim U(8,14)$;
$\ell_i\sim U(3.5,5)$; $e_i=s\,u_i$ with $u_i\sim U(0.5,1)$ and alternating signs.
The corridor is oriented so that the mean unconstrained-equilibrium flow is positive,
and the capacities are $60\%$ of its positive part (plus $0.05$); then the first-stage
constraint is active in all instances except $n=48$, and the second-stage one in
$40$--$100\%$ of the scenarios.

\subsection{Algorithms and implementation}\label{sec:algs}

All methods share the scenario decomposition, warm starts and the stopping criterion
$\mathrm{res}_k<10^{-3}$, where
$\mathrm{res}_k=\max_{\xi}\big|x_{k+1}(\xi)-\Pi_{C(\xi)}\big(x_{k+1}(\xi)-F(x_{k+1}(\xi),\xi)+w_{k+1}(\xi)\big)\big|_\infty$
is the natural residual of \eqref{eq:SVI}; since $x_{k+1}\in\NN$ and $w_{k+1}\in\MM$,
$\mathrm{res}_k=0$ if and only if $(x_{k+1},w_{k+1})$ solves \eqref{SVI2}.

\begin{itemize}
\item \textbf{IPHA} (baseline): progressive hedging \cite{RS:18}, i.e., $A=rI$ and
$x_{k+1}=P_\NN(\hat z)$, $w_{k+1}=w_k-rP_\MM(\hat z)$, with subproblems solved by
semismooth Newton. We report $r=1.1\max_\xi\norm{M_\xi}_2$, the value forced by
fixed-point inner solvers, and $r=r^*$, the best value in outer iterations on the grid
$r=f\max_\xi\norm{M_\xi}_2$, $f\in\{2^{-6},\dots,2\}\cup\{1.1\}$, selected a posteriori
for each instance (which favors the baseline). By Remark~\ref{rem:PDA}, this is AIPHA
with $A=rI$ and $S=P_\NN+r^2P_\MM$, so the comparison isolates the effect of the
nonsymmetric metric. Linear algebra is dense, except for IPHA with $r=r^*$ in P4.
\item \textbf{AIPHA}: Algorithm~\ref{AIPHA-eps} with $A=\Lambda-M_\xi$, $\Lambda$
diagonal, with $\Lambda_{ii}=1.05\sum_j|\sym(M_\xi)_{ij}|$ on the producer variables
(so $\sym(A)\succ0$ by Gershgorin's theorem, uniformly in the asymmetry) and
$\Lambda_{ii}=\norm{e}$ on the prices; in P3--P4 this is
$A_{\mathrm{red}}=\bigl(\begin{smallmatrix}\Lambda_x-G&-E\\E^\top&\norm{e}I\end{smallmatrix}\bigr)$.
The first-stage part of $\Lambda$ is scenario independent, so $A$ preserves $\NN$ and
$\MM$ (Remark~\ref{rem:twostage}). Since $M_\xi+A=\Lambda$, each subproblem is one
weighted projection onto $C(\xi)$: matrix-free and parallelizable across producers
and scenarios. We take $S=P_\NN+\kappa P_\MM$ with $\kappa=\norm{A}_2$ (better than
$\norm{A}_2^2$) and $\tau_k\equiv1$.
\item \textbf{Naive step}: the progressive-hedging step
$x_{k+1}=P_\NN(\hat z)$, $w_{k+1}=w_k-AP_\MM(\hat z)$ with the same $A$, which has no
convergence guarantee.
\end{itemize}
In Section~\ref{sec:results} the subproblems are solved exactly and $S$ is constant;
inexact subproblems are studied in Section~\ref{sec:inexact}. All instances, methods,
tables and figures are reproducible from the notebook \texttt{aipha\_pipeline.ipynb}.

\subsection{Results}\label{sec:results}

\begin{figure}[ht]\centering
\includegraphics[width=.95\textwidth]{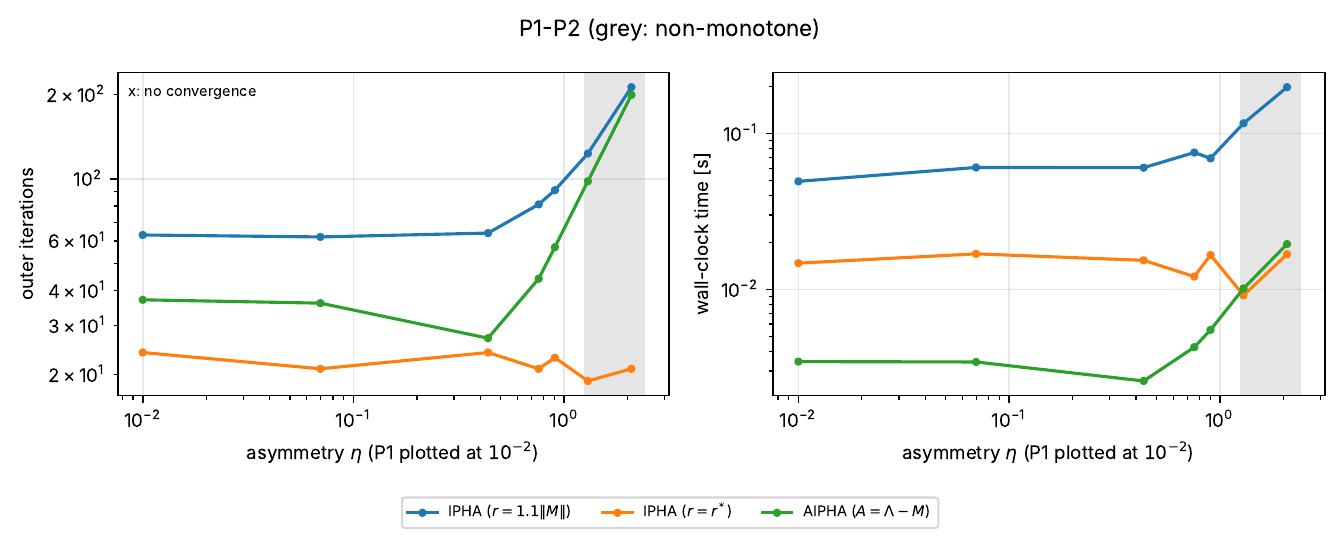}
\caption{P1--P2: outer iterations (left) and wall-clock time (right) vs.\ asymmetry $\eta$; the shaded region is not monotone.}\label{fig:P2}
\end{figure}

\paragraph{P1--P2.}
Table~\ref{tab:P2} and Figure~\ref{fig:P2} report P1 and the $\delta$-sweep of P2.
For progressive hedging, the choice of $r$ matters more than the asymmetry: with
$r^*\approx0.25\norm{M}$, IPHA needs $2.6$--$4$ times fewer outer iterations than with
$r=1.1\norm{M}$, and its iteration count ($21$--$24$) is essentially insensitive to
$\delta$. AIPHA is $2.8$--$5.7$ times faster than tuned IPHA on the monotone instances,
although it performs $1.1$--$2.5$ times more outer iterations, a number that grows with
$\eta$. On the non-monotone stress instances ($\delta\ge4$) this advantage disappears.

\begin{table}[ht]\centering\footnotesize\setlength{\tabcolsep}{2.3pt}
\caption{P1--P2: wall-clock time in seconds (outer iterations in parentheses), $\mathrm{res}_k<10^{-3}$. $r^*$ is the best progressive-hedging parameter on the grid $r=f\|M\|$. Best time in each column in bold. The columns $\delta\ge4$ are not monotone (stress test, outside the scope of the theory).}
\label{tab:P2}
\begin{tabular}{lrrrrrrr}
\toprule
 & P1 & $\delta=0$ & $\delta=1$ & $\delta=2$ & $\delta=2.5$ & $\delta=4$ & $\delta=8$\\
$\eta$ & 0.00 & 0.07 & 0.43 & 0.76 & 0.90 & 1.29 & 2.08\\
$r^*/\|M\|$ & 0.25 & 0.25 & 0.25 & 0.25 & 0.25 & 0.125 & 0.0625\\
\midrule
IPHA ($r=1.1\|M\|$) & 0.049\,(63) & 0.061\,(62) & 0.060\,(64) & 0.076\,(81) & 0.069\,(91) & 0.116\,(123) & 0.198\,(212)\\
IPHA ($r=r^*$) & 0.015\,(24) & 0.017\,(21) & 0.015\,(24) & 0.012\,(21) & 0.017\,(23) & \textbf{0.009}\,(19) & \textbf{0.017}\,(21)\\
AIPHA ($A=\Lambda-M$) & \textbf{0.003}\,(37) & \textbf{0.003}\,(36) & \textbf{0.003}\,(27) & \textbf{0.004}\,(44) & \textbf{0.006}\,(57) & 0.010\,(98) & 0.020\,(199)\\
\bottomrule
\end{tabular}
\end{table}

\begin{figure}[ht]\centering
\includegraphics[width=.95\textwidth]{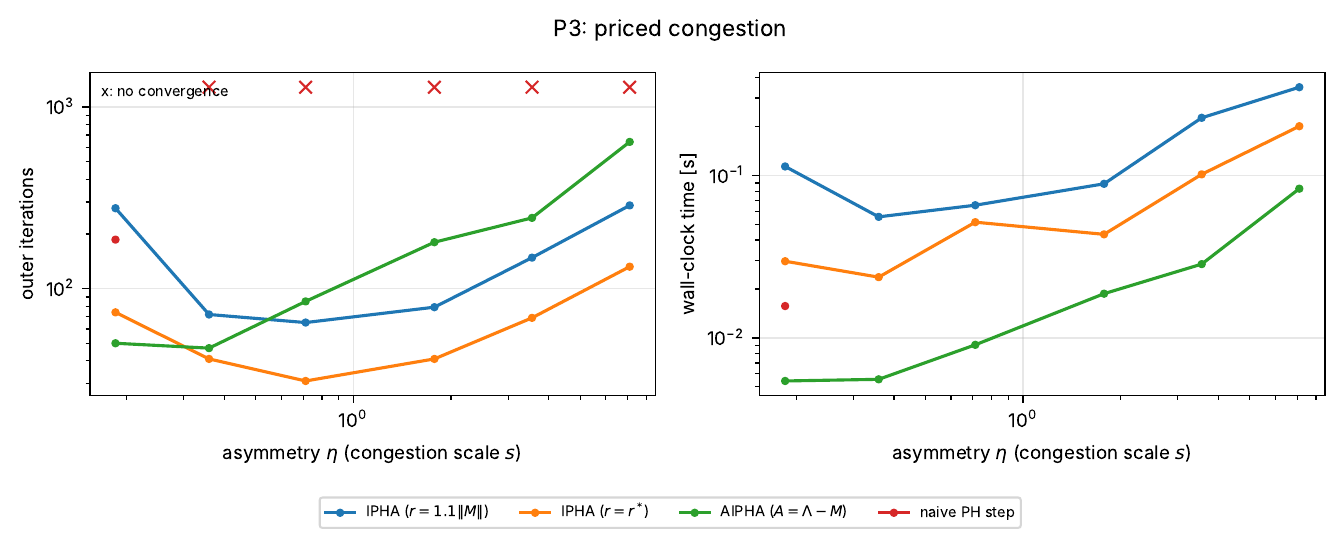}
\caption{P3: outer iterations (left) and wall-clock time (right) vs.\ asymmetry $\eta$ (congestion scale $s$); crosses mark runs that do not converge.}\label{fig:P3}
\end{figure}

\paragraph{P3: congestion.}
Table~\ref{tab:P3} and Figure~\ref{fig:P3} sweep the congestion scale
($\eta:0.18\to7.11$). AIPHA converges over the entire range and is the fastest method
for every $s$, $2.3$--$6.0$ times faster than tuned IPHA. Its iteration count grows
with $\eta$ (from $50$ to $642$, as $\norm{A_{\mathrm{red}}}\sim s$), whereas tuned
IPHA needs $31$--$132$, so the advantage shrinks as congestion grows. The naive step
with the same $A_{\mathrm{red}}$ fails to converge for every $s\ge1$ ($\eta\ge0.36$):
the outer step of Algorithm~\ref{AIPHA-eps} is \emph{not optional}, even in the
economically moderate regime.

\begin{table}[ht]\centering\footnotesize\setlength{\tabcolsep}{3.5pt}
\caption{P3: wall-clock time in seconds (outer iterations in parentheses) vs.\ congestion scale $s$; $n=4$. \textit{no conv.}\ means that the method did not reach the tolerance in 20000 iterations.}
\label{tab:P3}
\begin{tabular}{lrrrrrr}
\toprule
$s$ & 0.5 & 1 & 2 & 5 & 10 & 20\\
$\eta$ & 0.18 & 0.36 & 0.71 & 1.78 & 3.55 & 7.11\\
$r^*/\|M\|$ & 0.125 & 0.25 & 0.5 & 0.5 & 0.5 & 0.5\\
\midrule
IPHA ($r=1.1\|M\|$) & 0.114\,(277) & 0.056\,(72) & 0.066\,(65) & 0.089\,(79) & 0.227\,(148) & 0.350\,(287)\\
IPHA ($r=r^*$) & 0.030\,(74) & 0.024\,(41) & 0.052\,(31) & 0.043\,(41) & 0.102\,(69) & 0.201\,(132)\\
AIPHA ($A=\Lambda-M$) & \textbf{0.005}\,(50) & \textbf{0.006}\,(47) & \textbf{0.009}\,(85) & \textbf{0.019}\,(180) & \textbf{0.028}\,(245) & \textbf{0.083}\,(642)\\
naive PH step ($A_{\mathrm{red}}$) & 0.016\,(186) & no conv. & no conv. & no conv. & no conv. & no conv.\\
\bottomrule
\end{tabular}
\end{table}

\begin{figure}[ht]\centering
\includegraphics[width=.95\textwidth]{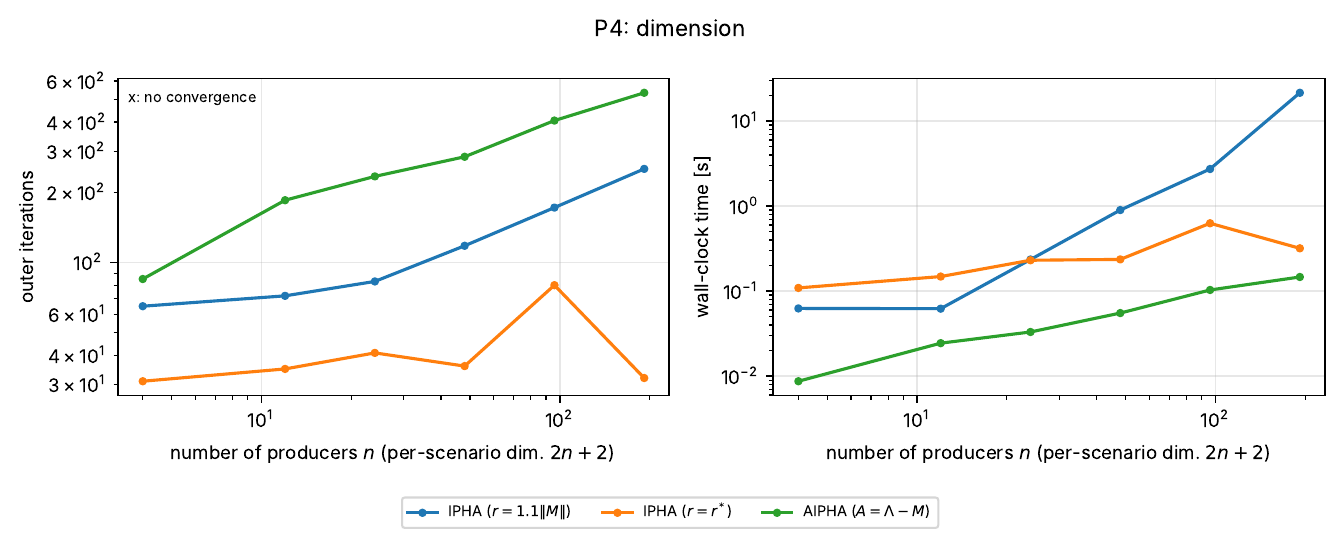}
\caption{P4: outer iterations (left) and wall-clock time (right) vs.\ number of producers $n$ (per-scenario dimension $2n+2$); the asymmetry $\eta$ grows with $n$.}\label{fig:P4}
\end{figure}

\paragraph{P4: dimension.}
Table~\ref{tab:P4} and Figure~\ref{fig:P4} grow the system ($\eta\sim\sqrt n$). With
dense Newton solves and $r=1.1\norm{M}$, IPHA grows from $0.06$ to $21.5$ seconds
($O(n^3)$ factorizations); exploiting the block-diagonal structure of $G_\xi$ with
sparse factorizations and $r=r^*$, it stays below $0.7$ seconds. AIPHA, with one
sparse matrix--vector product and one separable projection per iteration, remains the
fastest method for every $n$, $2.2$--$12.5$ times faster than tuned sparse IPHA despite
$2.7$--$17$ times more outer iterations; the advantage is not monotone in $n$ (at
$n=192$ tuned IPHA needs only $32$ iterations). When factorizations are unavailable or
the operator is only accessible through matrix--vector products, AIPHA needs nothing
beyond projections.

\begin{table}[ht]\centering\footnotesize\setlength{\tabcolsep}{3.5pt}
\caption{P4: nodal markets, wall-clock time in seconds (outer iterations in parentheses) vs.\ number of producers $n$. IPHA with $r=1.1\|M\|$ uses dense Newton solves; IPHA with $r=r^*$ uses sparse ones. The last row is the speedup of AIPHA over tuned sparse IPHA.}
\label{tab:P4}
\begin{tabular}{lrrrrrr}
\toprule
$n$ & 4 & 12 & 24 & 48 & 96 & 192\\
per-scenario dim.\ $2n{+}2$ & 10 & 26 & 50 & 98 & 194 & 386\\
$\eta$ & 0.71 & 1.13 & 1.72 & 2.28 & 3.15 & 4.43\\
$r^*/\|M\|$ & 0.5 & 0.5 & 0.5 & 0.25 & 0.5 & 0.125\\
\midrule
IPHA ($r=1.1\|M\|$) & 0.06\,(65) & 0.06\,(72) & 0.24\,(83) & 0.90\,(118) & 2.73\,(172) & 21.47\,(252)\\
IPHA ($r=r^*$) & 0.11\,(31) & 0.15\,(35) & 0.23\,(41) & 0.24\,(36) & 0.63\,(80) & 0.32\,(32)\\
AIPHA ($A=\Lambda-M$) & \textbf{0.01}\,(85) & \textbf{0.02}\,(185) & \textbf{0.03}\,(234) & \textbf{0.06}\,(284) & \textbf{0.10}\,(406) & \textbf{0.15}\,(534)\\
\midrule
speedup vs.\ IPHA ($r=r^*$) & $12.5\times$ & $6.1\times$ & $7.0\times$ & $4.3\times$ & $6.1\times$ & $2.2\times$\\
\bottomrule
\end{tabular}
\end{table}

\subsection{Inexact subproblems and number of scenarios}\label{sec:inexact}

To exercise the relative-error criteria, we run AIPHA with $A=r^*I$ and
$S=P_\NN+r^{*2}P_\MM$ (tuned progressive hedging for $\sigma=0$, by
Remark~\ref{rem:PDA}) and solve the subproblems by projected-gradient steps
$x\leftarrow\Pi_{C}(x-\beta(Hx-q))$. Given an approximate solution $\tilde x$, one more
step $\hat w=\Pi_C(u)$, $u=\tilde x-\beta(H\tilde x-q)$, gives
$n=(u-\hat w)/\beta\in N_C(\hat w)$, and $\hat x=x_k-(F(\hat w)+n-w_k)/r^*$ satisfies
\eqref{ips-ph} exactly with $\varepsilon_k=0$; inner steps are added until
\eqref{eq:conddeltaeps3}--\eqref{eq:conddeltaeps4} hold. With $\sigma\in[0.1,0.5]$
(Table~\ref{tab:inexact}) the inner steps drop by a factor of $12$--$17$ and the
wall-clock time by $6$--$8$, at the price of $5$--$80\%$ more outer iterations.
Finally, when the number of scenarios grows from $5$ to $200$ (P3 with $s=2$, P4 with
$n=24$), the outer iteration counts do not grow and the wall-clock times grow much
more slowly than the number of scenarios (vectorized subproblems); with $200$
scenarios AIPHA is $5$--$7$ times faster than tuned IPHA.

\begin{table}[ht]\centering\small\setlength{\tabcolsep}{4pt}
\caption{Inexact AIPHA with $A=r^*I$, $S=P_{\mathcal N}+r^{*2}P_{\mathcal M}$: outer iterations, total inner projected-gradient steps and wall-clock time (s) vs.\ the relative-error tolerance $\sigma$ ($\sigma=0$: exact subproblems, i.e., tuned progressive hedging).}
\label{tab:inexact}
\begin{tabular}{crrrrrrrrr}
\toprule
 & \multicolumn{3}{c}{P3, $s=2$} & \multicolumn{3}{c}{P3, $s=10$} & \multicolumn{3}{c}{P4, $n=24$}\\
\cmidrule(lr){2-4}\cmidrule(lr){5-7}\cmidrule(lr){8-10}
$\sigma$ & outer & inner & time & outer & inner & time & outer & inner & time\\
\midrule
0 & 31 & 3836 & 0.201 & 69 & 7634 & 0.394 & 41 & 4900 & 0.377\\
0.1 & 33 & 326 & 0.035 & 75 & 646 & 0.070 & 43 & 416 & 0.058\\
0.3 & 40 & 252 & 0.029 & 89 & 454 & 0.068 & 51 & 300 & 0.045\\
0.5 & 52 & 232 & 0.027 & 116 & 444 & 0.057 & 73 & 286 & 0.049\\
0.7 & 70 & 222 & 0.028 & 171 & 422 & 0.072 & 105 & 270 & 0.051\\
0.9 & 110 & 220 & 0.030 & 212 & 424 & 0.092 & 134 & 268 & 0.053\\
\bottomrule
\end{tabular}
\end{table}

\subsection{Discussion}\label{sec:discussion}

With a scalar metric, AIPHA is progressive hedging, and the relevant baseline uses a
well-chosen $r$: the value $r\approx\norm{M}$ imposed by fixed-point inner solvers is
far from optimal. The defining freedom of AIPHA is that $A$ may be
\emph{nonsymmetric}: only a nonsymmetric $A$ can cancel the skew part of $M$ inside the
subproblem, since a symmetric $A$ always leaves $\skw(M)$ in the resolvent. With
$A=\Lambda-M$, every subproblem becomes a closed-form, matrix-free, decentralized
projection. The price is an outer iteration count that grows with $\norm{A}$, but each
iteration is so cheap that AIPHA is faster than tuned progressive hedging on all
monotone instances, provided the safeguarded outer update of Algorithm~\ref{AIPHA-eps}
is used.

\section{Conclusions}\label{sec:conclusions}

We have introduced AIPDA, an asymmetric inexact progressive decoupling algorithm for
linkage monotone inclusions, and shown that it is an exact instance of the asymmetric
variable-metric proximal point framework of \cite{BLPS:25} applied to the partial
inverse; it therefore inherits global convergence and a local linear rate, the latter
without additional assumptions when the operator is polyhedral. AIPDA recovers
classical progressive decoupling and progressive hedging with scalar metrics, contains
the inexact progressive decoupling algorithm of \cite{rock:24,Rock:26} from some
iteration on under a blockwise error condition, and extends that family with
nonsymmetric metrics, a free parameter $c_k$ and relative errors. Its stochastic
specialization AIPHA decomposes by scenarios; for two-stage affine problems without
first-stage/recourse coupling, the nonsymmetric metric $A=\Lambda-M$ reduces every
scenario subproblem to a single weighted projection. On stochastic Cournot electricity
markets with priced congestion, this makes AIPHA faster than tuned progressive hedging
on all monotone instances, and the relative-error criteria reduce the inner work by an
order of magnitude.

\backmatter

%
\section*{Declarations}

\begin{itemize}
\item \textbf{Funding.} This work was partially supported by  Universidad Nacional de Rosario
(projects~80020220700214UR and 80020240200045UR).
\item \textbf{Competing interests.} The authors have no competing interests
to declare that are relevant to the content of this article.
\item \textbf{Data availability.} No external datasets were used. All problem
instances in Section~\ref{sec:numerics} are synthetic and are generated
reproducibly (fixed random seed) by the code described below.
\item \textbf{Code availability.} The code implementing the algorithms and
reproducing all experiments, tables and figures (Jupyter notebook
\texttt{aipha\_pipeline.ipynb}) is available from the corresponding author
upon reasonable request.
\end{itemize}

\bibliography{AIPDA}

\end{document}